 \documentclass[11pt,a4paper]{article}
\usepackage[percent]{overpic}
\usepackage[T1]{fontenc}
\usepackage{wrapfig,pgfplots,caption,sidecap,subcaption,wrapfig,graphicx,lipsum,hyperref,titlesec,verbatim,indentfirst,amsmath, amsthm, amssymb, textcomp, enumerate, multicol, fancyhdr, varwidth, graphicx, color, mathrsfs, mathtools, pdflscape, accents, mdframed,soul, mathptmx, relsize, tikz, float, nonfloat, amsfonts, sectsty,graphicx, euscript, setspace,ar, lineno, blindtext,lmodern}
\DeclareMathAlphabet{\mathpzc}{T1}{pzc}{m}{it}
\SetMathAlphabet\EuScript{bold}{U}{eus}{b}{n}
\usepackage[ruled,vlined]{algorithm2e}
\usepackage{todonotes}
\newcommand{\etal}{\textit{et al.} }
\newcommand{\nmesh}{$\widehat{\mathbf{M}}(t)$ }
\usepackage[font=small,labelfont=bf]{caption}
\usepackage[export]{adjustbox}
 \usepackage[left=25mm,right=25mm,top=28.5mm,bottom=28.5mm]{geometry}
\usepackage[outline]{contour}
\title{A simple history dependent remeshing technique to increase finite element model stability in elastic surface deformations}
\author{Jessica R. Crawshaw$^1$, Jennifer A. Flegg$^1$, James M. Osborne$^1$\\
$^1$School of Mathematics and Statistics, The University of Melbourne, Melbourne, Australia}
\graphicspath{{Figures/} }
\date{}  
 \usepackage[numbers,sort&compress]{natbib}
\definecolor{purple}{rgb}{0.5, 0, 0.5}
\definecolor{green}{rgb}{0, 0.7, 0}
\definecolor{blue}{rgb}{0.2, 0.3, 1}
\definecolor{orange}{rgb}{0.95, 0.61, 0}
  \definecolor{lblue}{rgb}{0.2, 0.3, 1}
\def\*#1{\mathbf{#1}}
 \DeclareMathAlphabet      {\mathbfit}{OML}{cmm}{b}{it}
\renewcommand{\vec}[1]{\mathbfit{#1} } 
 \usepackage[noend]{algpseudocode}
\begin{document}
 \maketitle
\begin{abstract}
In this paper, we present and validate a simple adaptive surface remeshing technique to transfer history dependent variables from an old distorting mesh to a new mesh during finite element simulations of elastic surface deformation. This technique allows us to reduce the error arising from excessive mesh distortion whilst preserving information about the initial configuration of the mesh and the history dependent variables. The transfer technique presented here constructs the initial configuration of the new mesh by considering the distortion incurred by the elements of the old mesh and projecting backward in time. Using this new initial configuration, the stress and strain over the new mesh can be easily calculated. 
After presenting the necessary steps to reconstruct the initial configuration, we show that this relatively simple transfer technique adds stability to finite element simulations and reduces the spatial error and the strain error across the domain. The novel transfer technique presented in this paper is easy to implement and provides a simple strategy to add stability to simulations undergoing large deformations.
\end{abstract}

	\section{Introduction}
\begin{figure}[t]
	\centering
			\begin{subfigure}[b]{0.67\textwidth}
					\hspace{-1.2cm}
					 \vspace{-0.5cm}
					\centering
					\begin{overpic}[width=\textwidth]{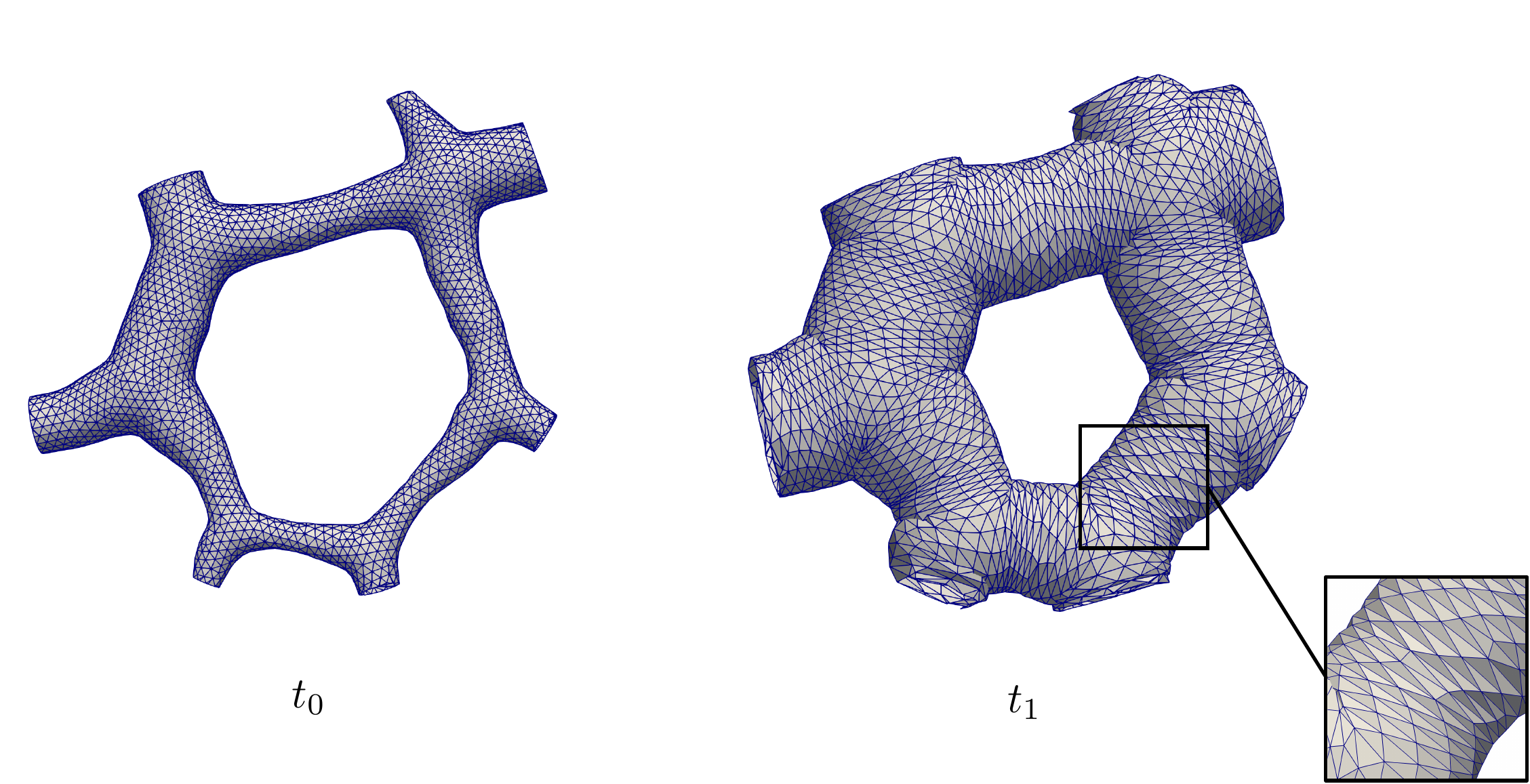}
					\put (04, 49) {$a)$}
					\put (50, 49) {$b)$}
					\put (103, 49) {$c)$}
					\put(19, 5){\color{white}\circle*{10}}
					\put(19, 5){$t_0$}
					\put(65.5, 5){\color{white}\circle*{10}}
					\put(65.5, 5){$t_1$}
				\end{overpic}
			\end{subfigure}
			 \vspace{1cm}
			\begin{subfigure}[b]{0.3\textwidth}
					\centering
					\begin{overpic}[width=0.9\textwidth, trim={1.6cm 3.7cm 0cm 5cm},clip, right]{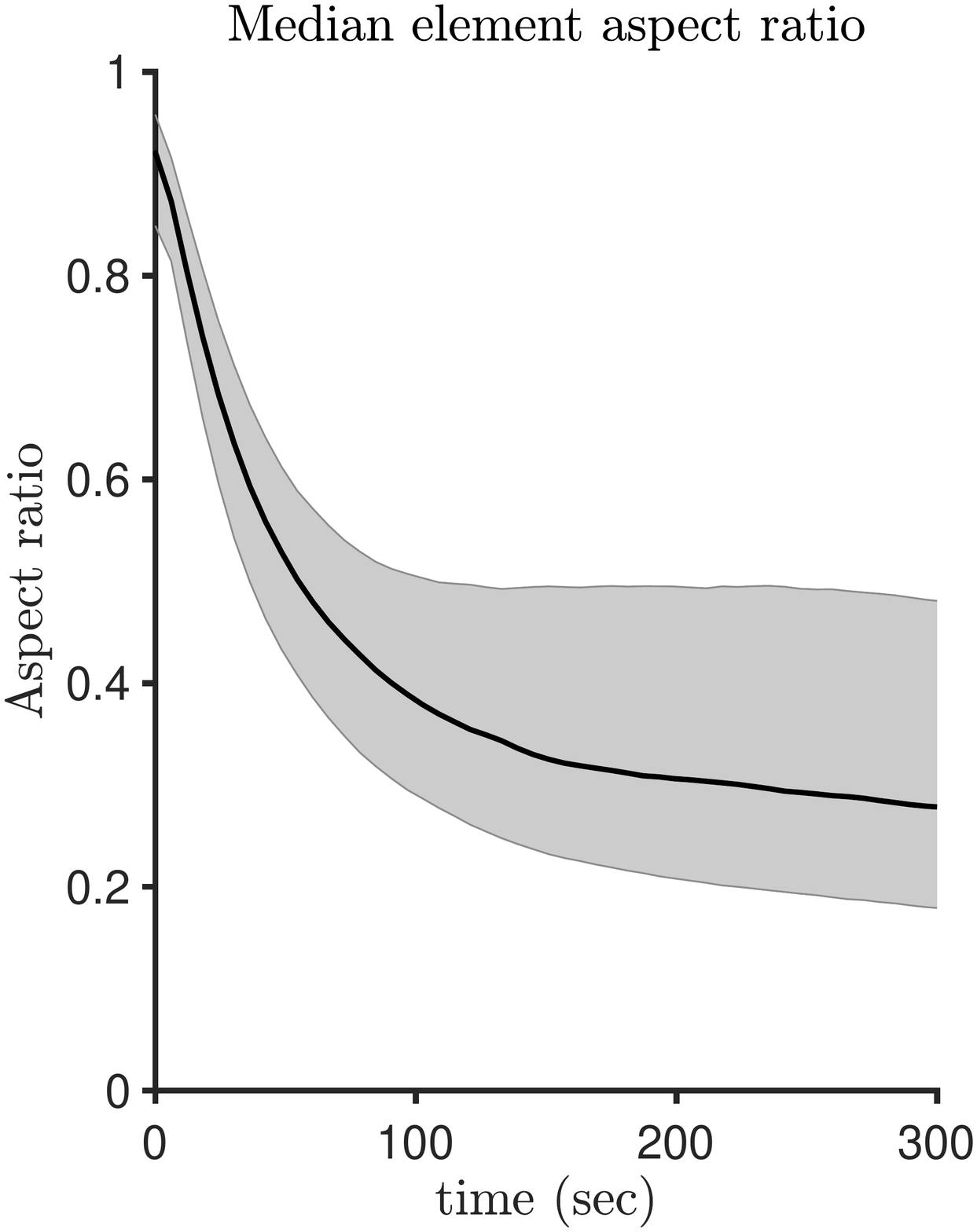}
					 \put (24, 103) {\fontsize{9}{9}\selectfont Element aspect ratio}
					 \put (41, 95) {\fontsize{9}{9}\selectfont over time}
					 \put(42,-9) {\fontsize{9}{9}\selectfont time (sec)}
					\put(1.7, 34){\rotatebox{90}{ \fontsize{9}{9}\selectfont aspect ratio}}
					\end{overpic}
			\end{subfigure}
			\caption{\setstretch{1.4}\textbf{A finite element model of vascular deformation with a fixed mesh showing mesh degradation.} A constant internal pressure deforms a capillary plexus from an initial configuration at $t_0$ ($a$) to a deformed state at time $t_1$ ($b)$. As the plexus grows, the mesh progressively coarsens, and the elements become increasingly distorted. $c)$ The median element aspect ratio (solid black line) drops over the duration of the simulation, while the interquartile range (shaded region) expands.}\label{fig:ExampleOfMeshDegredation}
	\end{figure}

	 When investigating large geometric deformations using finite element models it is critical to maintain mesh integrity to ensure model accuracy \cite{luo2008dealing}.  Accuracy depends on both the refinement and the regularity of the mesh used to discretise the deforming domain. 
	 This is particularly important during large deformations in which increasingly excessive mesh coarsening and element distortion degrades the quality of the mesh throughout the simulation \cite{javani2014consistent}. 
	 As an example, Figure \ref{fig:ExampleOfMeshDegredation}a shows a deforming capillary plexus, discretised using an unstructured triangular mesh (the plexus is derived from images in Ghaffari \etal \cite{ghaffari2015simultaneous} and the deformation is based on the model presented in Kr{\"u}ger \etal \cite{kruger2011efficient}, using the algorithm from Osborne and Bernabeu \cite{osborne2018fully}). As the pressurised plexus grows, the mesh becomes increasingly distorted and elements elongate in the tangential direction (Figure \ref{fig:ExampleOfMeshDegredation}b). Excessive element distortion over a mesh, as shown in Figure \ref{fig:ExampleOfMeshDegredation}c, introduces large, unacceptable error to the simulation, thereby jeopardising the accuracy of the model and possibly invalidating any further analysis.
	 To overcome these complications and to continue the simulation, successive remeshing throughout the simulation is unavoidable \cite{zeramdini2019numerical}.\\

	 Following the generation of a new mesh in the remeshing process, the simulation variables from the old degenerated mesh must be transferred to the new mesh. One can either completely recompute the simulation variables or transfer the variables from the old mesh to the new mesh \cite{zeramdini2019numerical}. If the optimal mesh configuration changes continuously thoughout the simulation (as is the case in the example shown in Figure \ref{fig:ExampleOfMeshDegredation}) the second approach is generally considered preferable. However, doing so is a delicate matter, as inadequate determination of the updated state variables can introduce remeshing error which will adversely affect the simulation \cite{dureisseix2006information}. This process of transferring the history-dependent variables from the old to new mesh is known as history-dependent remeshing. History-dependent remeshing techniques have been developed over the past 30 years prompted by the need for increasing model accuracy \cite{dureisseix2006information,srikanth2000shape,  kumar2015parallel,hinton1974local,peric1996transfer, zeramdini2019numerical, zienkiewicz1992superconvergent, zienkiewicz1992superconvergent2, zienkiewicz1992superconvergent3}. These transferal techniques can be broadly categorised into two groups; (1) techniques to transfer discontinuous variables which are continuous within each element and stored at Gauss points (referred to here as $P_0$ techniques) such as stress and strain ($P_0$ variables); and (2) techniques to transfer continuous variables stored at nodes  (referred to here as $P_1$ techniques) such as displacement, velocity and temperate ($P_1$ variables) \cite{kumar2015parallel}.\\

 $P_1$ transfer techniques are more straight forward than $P_0$ techniques, and are appropriate when describing continuous fields where information is stored at the nodes. These techniques typically utilise interpolation or extrapolation techniques to attain a continuous field over the original mesh from the nodal values \cite{dureisseix2006information}. This continuous field can then be used to compute the nodal values over the new mesh. Various papers have extended upon this simple concept to increase accuracy and robustness of the technique using the least squares method and shape functions \cite{bussetta2010comparison, kumar2015parallel, pere2002mapping}.\\

Owing to the discontinuous nature of $P_0$ variables such as stress and strain, $P_1$ techniques are inappropriate to transfer $P_0$ variables between meshes. Indeed, the transfer of $P_0$ variables requires considerably more sophisticated techniques. 
$P_0$ techniques typically (i) enforce the balance equation in a weak or strong sense, (ii) use first or second order interpolation techniques, or (iii) are based on nodal or element patches \cite{kumar2015parallel}. Following the first approach (i), Javani \etal \cite{javani2014consistent} transferred a minimum set of variables from the old to the new mesh, and then reconstructed the remaining variables via the governing constitutive and balance equations. This procedure ensured the consistency of the constitutive and equilibrium equations were preserved over the new mesh.\\

 First and Second order $P_0$ interpolation techniques (ii) generally follow one of two approaches; (a) direct interpolation \cite{kucharik2003efficient, liszka1984interpolation}; and (b) nodal recovery based techniques \cite{peric1996transfer, lee1994error, dureisseix2006information}.
Direct interpolation (a) approximates piecewise continuous fields over the new mesh by using a set of field values at Gauss points on the old mesh in the region of the considered Gauss point of the new mesh to build a local and continuous description of the variable \cite{lobov2011dll4,liszka1984interpolation}. An example of the direct interpolation method was presented by Brancherie \etal \cite{brancherie2008consistent}, who proposed a field transfer operator based the Moving Least Squares and the diffuse approximation methods to reconstruct the stress field over the entire domain. Nodal recovery based techniques (b) transfer the $P_0$ variables by firstly projecting the values of the Gauss points to the nodal points of the old mesh and build a continuous interpolation of the $P_0$ variable. From this continuous interpolation, the value at the nodes of the new mesh can be calculated. The $P_0$ variables at the Gauss points of the new mesh are then obtained from the nodal points using shape functions. Nodal recovery based techniques may not always be appropriate; in 2007 Khoei \etal \cite{khoei2007superconvergence} found that these techniques violate the equilibrium of the systems, even if the mesh is unchanged.\\

One of the most popular $P_0$ transfer techniques is the superconvergent patch recovery (SPR) technique (iii) proposed by Zienkiewicz \etal \cite{zienkiewicz1992superconvergent,zienkiewicz1992superconvergent2} in 1992. Using the SPR technique, for each node, a continuous polynomial expansion is taken over a patch of elements sharing the given node to construct a global approximation of the variables across the mesh. At points with superconvergent properties, this polynomial expansion approaches the values provided by the finite element analysis. An example of this technique in practise was presented by Boussetta \etal \cite{boussetta2006adaptive} who used the SPR technique to produce a continuous description of the stress tensor over the domain, from which the required mechanical properties could be reconstructed across the new mesh. 
Over the years, many modified and alternative SPR techniques have been proposed to increase the accuracy and robustness of this technique. Boroomand \etal \cite{boroomand1997recovery} proposed the Recovery by Equilibrium in Patch  technique. This technique is a  modification of the classic SPR technique that avoids the need to specify superconvergent points. Gu \etal \cite{gu2004modified} improved the performance of the SPR technique for nonlinear problems by using integration points as sampling points, rather than the traditional use of nodes as sampling points. Wiber \etal \cite{wiberg1995improved, wiberg1997superconvergent} formulated the superconvergent patches by elements transfer technique by introducing element patches rather than nodal patches. This technique penalised violation of the equilibrium and included a least squares fit for the boundary conditions. \\

It is evident that $P_1$ techniques are far more simple and straight forward than $P_0$ techniques. Whilst robust, the currently available $P_0$ techniques are often complicated, computationally time consuming, or involve a significant implementation time \cite{peric1996transfer,zienkiewicz1992superconvergent, hinton1974local, zienkiewicz1992superconvergent2, zienkiewicz1992superconvergent3}. 
  Unfortunately stress and strain are defined at each element meaning $P_1$ techniques are ill-equipped to deal with the discontinuous nature of these variables, thus forcing the use of $P_0$ techniques. 
In this paper we present a simple alternative to the classic $P_0$ and $P_1$ techniques to transfer strain from a deforming mesh to the new mesh in finite element models. 
This is done by focusing our attention on the displacement of each node and the continuous distortion incurred by each element, $P_1$ variables, over the old mesh. From this we can construct a new initial configuration for the new mesh over the original geometry, $\Omega_0$, and then calculate the strain over the new mesh (a $P_0$ variable) without resorting to more complicated $P_0$ techniques. This novel history-depended remeshing technique is of particular use in finite element models of elastic deformation where the reverse deformation is not known. We have designed this novel technique to be relatively straight forward to implement and appropriate for iterative remeshing schemes. In the remainder of this paper we describe this novel technique in more detail, examine its accuracy, and demonstrate how it can be implemented in a model of vascular deformation.
\section{Computational methodology} \label{Methods}
		\begin{figure}[b!]
				\centering
		 \vspace{0.5cm}
		 \begin{overpic}[width=0.6\textwidth, trim={5cm 6cm 1cm 4.5cm},clip]{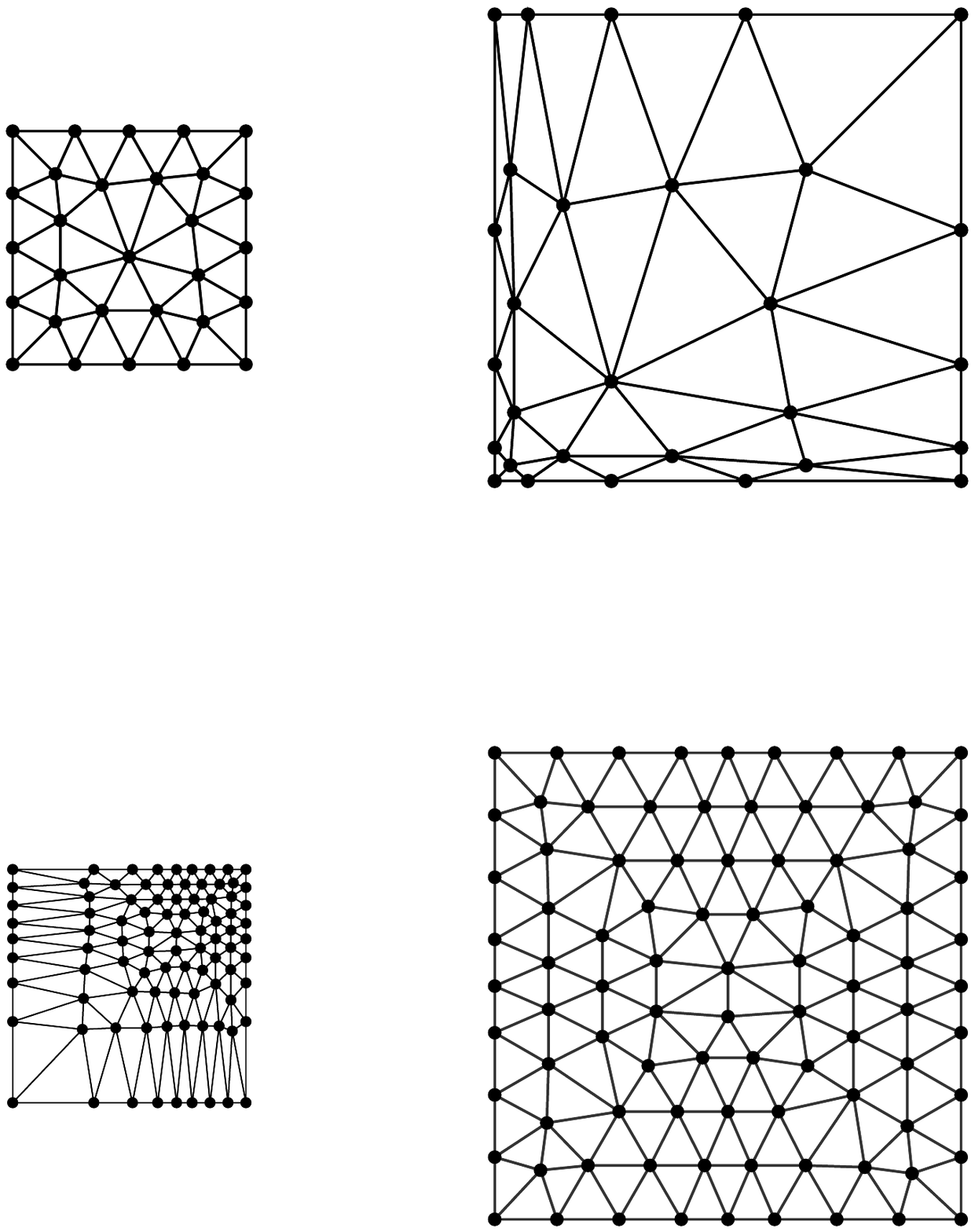}
				\put (7, 87) {$\mathbf{M}_0$}
				\put (51, 96) {$\mathbf{M}$} 
				\put (8, 32.5) {$\widehat{\pmb{\EuScript{M}}}_0$} 
				\put (51, 42) {$\widehat{\mathbf{M}}$}
				 \put (24, 22) {$\boldsymbol{\chi}(\widehat{\vec{X}}_i)$}
				 \put (25, 18) {$\mathrel{\contour{black}{${\longleftarrow}$}}$}
				 \put (24, 77) {$\vec{f}(\vec{x}_i)$}
				 \put (25, 73) {$\mathrel{\contour{black}{${\longrightarrow}$}}$}
				 \put (52, 50) {$\mathrel{\contour{black}{${\downarrow}$}}$}
				\end{overpic}
		 \caption{\setstretch{1.4} \textbf{An example of the history dependent variable transfer technique following remeshing of a 2D deforming square.} The initial geometry, $\Omega_0$, is discretised by $\mathbf{M}_0$ (top left). The square is deformed by $\vec{f}(\vec{x}_i)$, resulting in a non-uniform mesh with distorted elements, $\mathbf{M}$ (top right). In this example $\vec{f}(\vec{x}_i)= [x_i^2,y_i^2]$, however the deformation is usually unknown. Here $\vec{x}_i = [x_i,y_i]$ is the position of Node $i$ in the initial configuration and $\vec{X}_i = [X_i,Y_i]$ is the position of Node $i$ in the deformed mesh. The deformed square, given by the boundary of $\mathbf{M}$, is remeshed with a new, independent mesh, $\widehat{\mathbf{M}}$ (bottom right). History-dependent simulations require a knowledge of the initial configuration of the mesh. As such the initial positions for each node in $\widehat{\mathbf{M}}$ are attained via a mapping, $\boldsymbol{\chi}(\widehat{\vec{X}}_i)$, of $\widehat{\mathbf{M}}$ to the initial geometry producing $\widehat{\pmb{\EuScript{M}}}_0$ (bottom left), where $\widehat{\vec{X}}_i$ is the location of Node $i$ in the new mesh.}\label{fig:OldAndNewMesh}
	\end{figure}

	We begin with some initial geometry, $\Omega_0$, discretised by an initial mesh, $\mathbf{M}_0$. This mesh is   comprised of a set of nodes, $\mathbf{N}_0$, with locations, $\vec{x}_i = [x_i,y_i]$,  and elements, $\mathbf{E}$, (Figure \ref{fig:OldAndNewMesh}, $\mathbf{M}_0$). The mesh is subject to a deformation, $\vec{f}(\vec{x}_i)$, which is usually unknown, leading to a new configuration at a later time $t$ \footnote{  Note that we don't explicitly consider time $t$ at this point. Here we are considering only the initial state and the deformed state. Later we consider deformations over time.}, $\mathbf{M}$, with node locations $\vec{X}_i = [X_i,Y_i]$ (Figure \ref{fig:OldAndNewMesh},  $\mathbf{M}$). The elements in $\mathbf{M}$ are distorted and no longer provide an appropriate approximation of the deformed geometry.
	This deformed geometry is then remeshed with a completely new set of nodes and elements to provide a more accurate representation of the geometry (Figure \ref{fig:OldAndNewMesh}, $\widehat{\mathbf{M}}$).
	 The new mesh is denoted as $\widehat{\mathbf{M}}$, and is comprised of a new set of nodes and elements which are given by $\widehat{\mathbf{N}}$ and $\mathbf{\widehat{E}}$, respectively where $\widehat{\vec{X}}_i = [\widehat{X}_i,\widehat{Y}_i]$ are the new node locations. In history-dependent models, such as models of elasticity, the state variables across the domain are dependent on a clearly defined initial configuration. As such, we must ensure that each node in the new mesh, $\widehat{\mathbf{N}}$, has an initial position defined on the undeformed geometry, $\Omega_0$, which we denote as $\widehat{\pmb{\EuScript{N}}}_0$ where the initial node locations are given by $\widehat{\vec{x}}_i = [\widehat{x}_i,\widehat{y}_i]$. The set of nodes $\widehat{\pmb{\EuScript{N}}}_0$ and elements $\widehat{E}$ define the undeformed new mesh $\widehat{\pmb{\EuScript{M}}}_0$  (Figure \ref{fig:OldAndNewMesh}, $\widehat{\pmb{\EuScript{M}}}_0$).
	\begin{figure}[htb!]
				\centering
				\vspace{0.5cm}
				\begin{overpic}[width=0.5\textwidth]{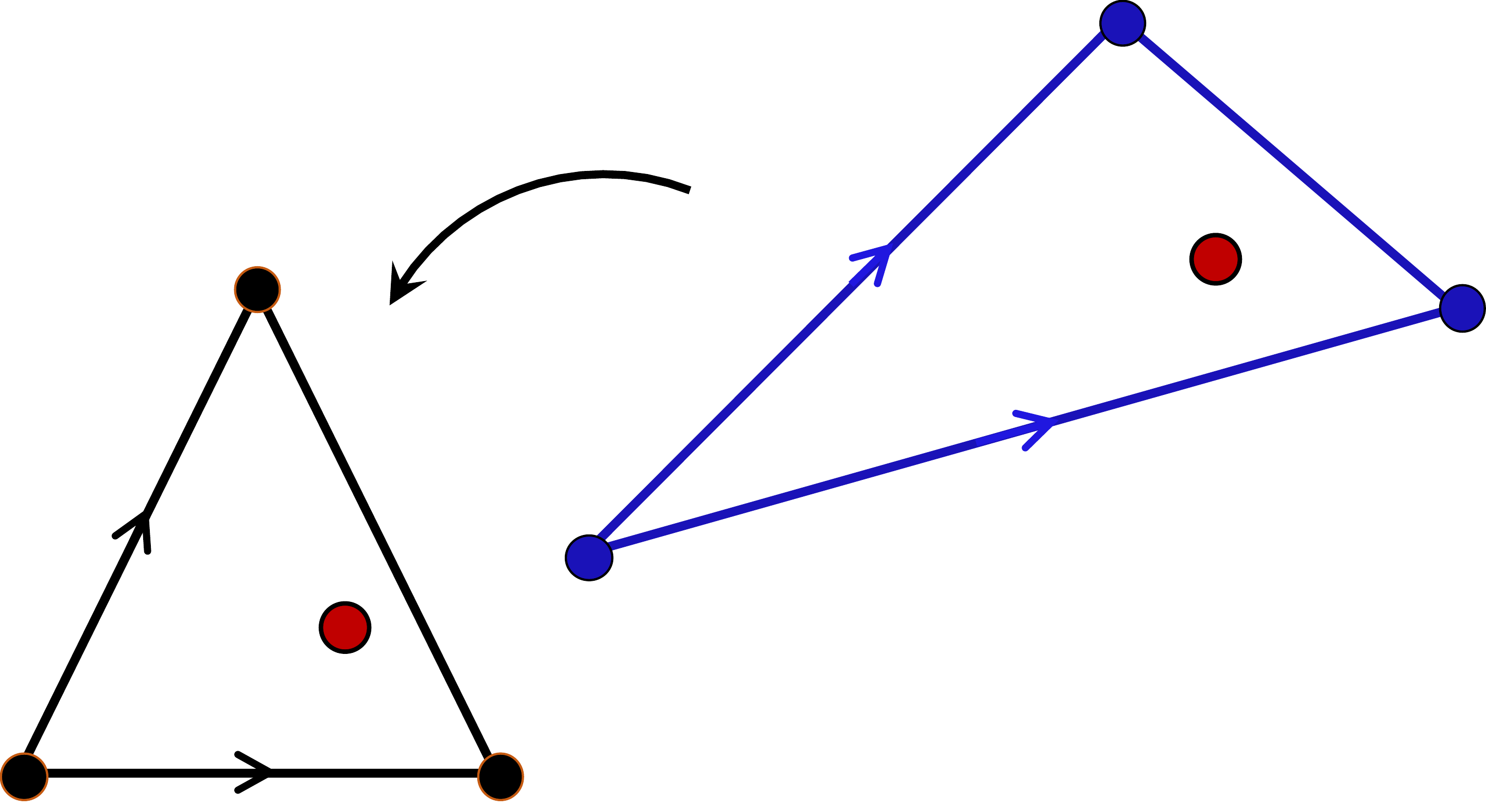}
				\put (-5, 18) {$\vec{u}^j_{0,2}$}
				\put (12.5, -5.5) {$\vec{u}^j_{0,1}$}
				\put (52, 42) {$\vec{u}^j_2$}
				\put (70, 20) {$\vec{u}^j_1$}
				\put (74, 40) {$\widehat{\vec{X}}_i$}
				\put (15 ,12) {$\widehat{\vec{x}}_i$}
				\put (1, -14.5) {Original element}
				\put (54, 8) {Deformed element}
				\end{overpic}
				\vspace{1.5cm}
					\caption{\setstretch{1.4}\textbf{A schematic of the transfer process for a single node.} Node $i$ (red) from the new mesh, located at $\widehat{\vec{X}}_i$, is shown in the deformed configuration within the nearest element, $j$, from the old mesh (blue element, right). The relative position of Node $i$ within Element $j$ is expressed using a change of basis with the edge vectors $\vec{u}_1^j$ and $\vec{u}_2^j$ and the element normal (not shown, Equation (\ref{eq:ChangeOfBasis})). By replacing these basis vectors with the equivalent vectors ($\vec{u}_{0,1}^j$ and $\vec{u}_{0,2}^j$) from the undeformed state (black element, left), the position of Node $i$, $\widehat{\vec{x}}_i$, over the original configuration can be determined.}\label{fig:SampleOfRemeshing}
				\end{figure}
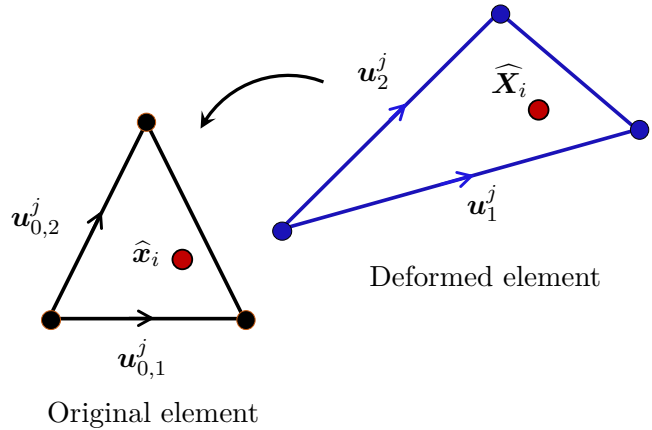
				
	\subsection{Mapping the new mesh to the  undeformed configuration}
	To determine $\widehat{\pmb{\EuScript{N}}}_0$, we map each node in $\widehat{\mathbf{N}}$ to the undeformed geometry, $\Omega_0$. This is achieved by considering the relative position of each node within the nearest element, $j$, in the old mesh, $\mathbf{M}$, and the deformation incurred by this element (Figure \ref{fig:SampleOfRemeshing}). The method to calculate the  nearest element is discussed below (Section \ref{NearestElementSearch}). It is assumed that the deformation within each element is linear \cite{kruger2011efficient,newham2019finite,zauel2006comparison}, however this assumption could be relaxed with further model development (see Discussion). The relative position, $\widehat{\vec{X}}_i$, of Node $i$ in $\widehat{\mathbf{M}}$  is expressed in the basis defined by three linearly independent vectors characterising Element $j$,
		\begin{equation}\label{eq:ChangeOfBasis}
				\widehat{\vec{X}}_i= c_1\vec{u}^j_1+ c_2\vec{u}^j_2+ c_3\vec{n}^{j}.
		\end{equation}
		The three vectors, $\vec{u}^j_1$, $\vec{u}^j_2$ and $\vec{n^j}$, uniquely characterise Element $j$ (Figure \ref{fig:SampleOfRemeshing}) and represent two edge vectors and the unit normal at time $t$, respectively. Here the unit normal is defined as $\vec{n^j}= \vec{u}^j_1 \times \vec{u}^j_2/|\vec{u}^j_1 \times \vec{u}^j_2|$. The basis coefficients, $c_1,~c_2$ and $c_3$ are determined from the set of simultaneous equations arising from a basis transformation from Cartesian coordinates to the new basis defined by  Element $j$. As we assume deformation across each element is linear, the deformation at each point within each element at any time can be described using the deformation of the basis vectors. Thus it follows that we can describe the position of Node $i$ in $\Omega_0$ by interchanging the basis vectors in Equation (\ref{eq:ChangeOfBasis}) with the basis vectors characterising Element $j$ in the undeformed state,
		\begin{equation}  
		\widehat{\vec{x}}_i= c_1\vec{u}^j_{0,1} + c_2\vec{u}^j_{0,2}+ c_3\vec{n}^j_0.
		\end{equation}
		Here  $\vec{u}^j_{0,1}$, and $\vec{u}^j_{0,2}$ are two edge vectors of Element $j$ in the undeformed state, and $\vec{n^j}_0$ is the unit normal in this state.  
		Using this relatively simple method, we are able to map each node from the new mesh, $\widehat{\mathbf{M}}$, to the initial geometry, $\Omega_0$, thus providing an initial configuration for each new element. This mapped initial configuration allows us to recalculate the history-dependent state variables, such as stress and strain, over the new mesh. In doing so, we are able to remesh the history-dependent domain during large elastic deformation and maintain a near optimal discretisation of the domain.  

\subsection{Nearest element search algorithm} \label{NearestElementSearch}
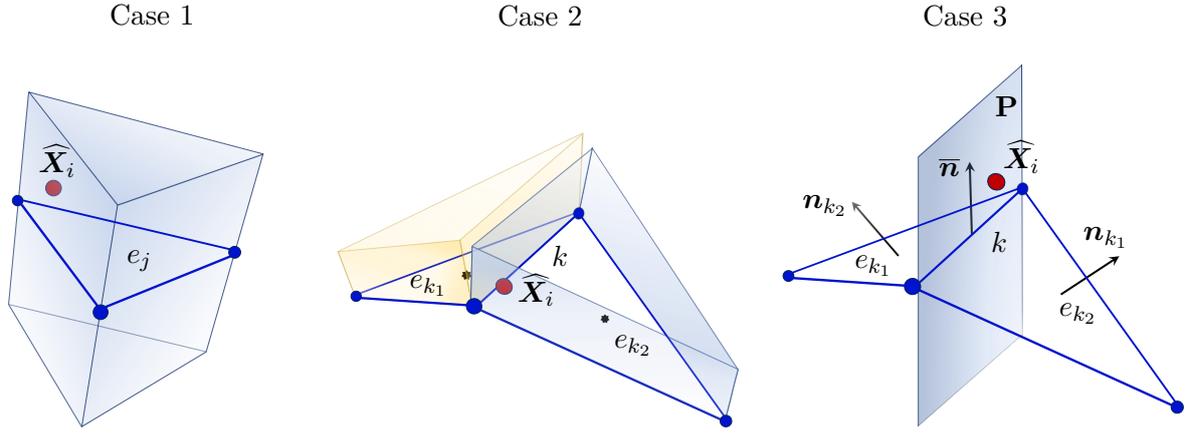
\begin{figure}[t!]
			\centering
			\begin{subfigure}[b]{0.3\textwidth}
			 \hspace{-0.4cm}
			\begin{overpic}[width=0.7\textwidth]{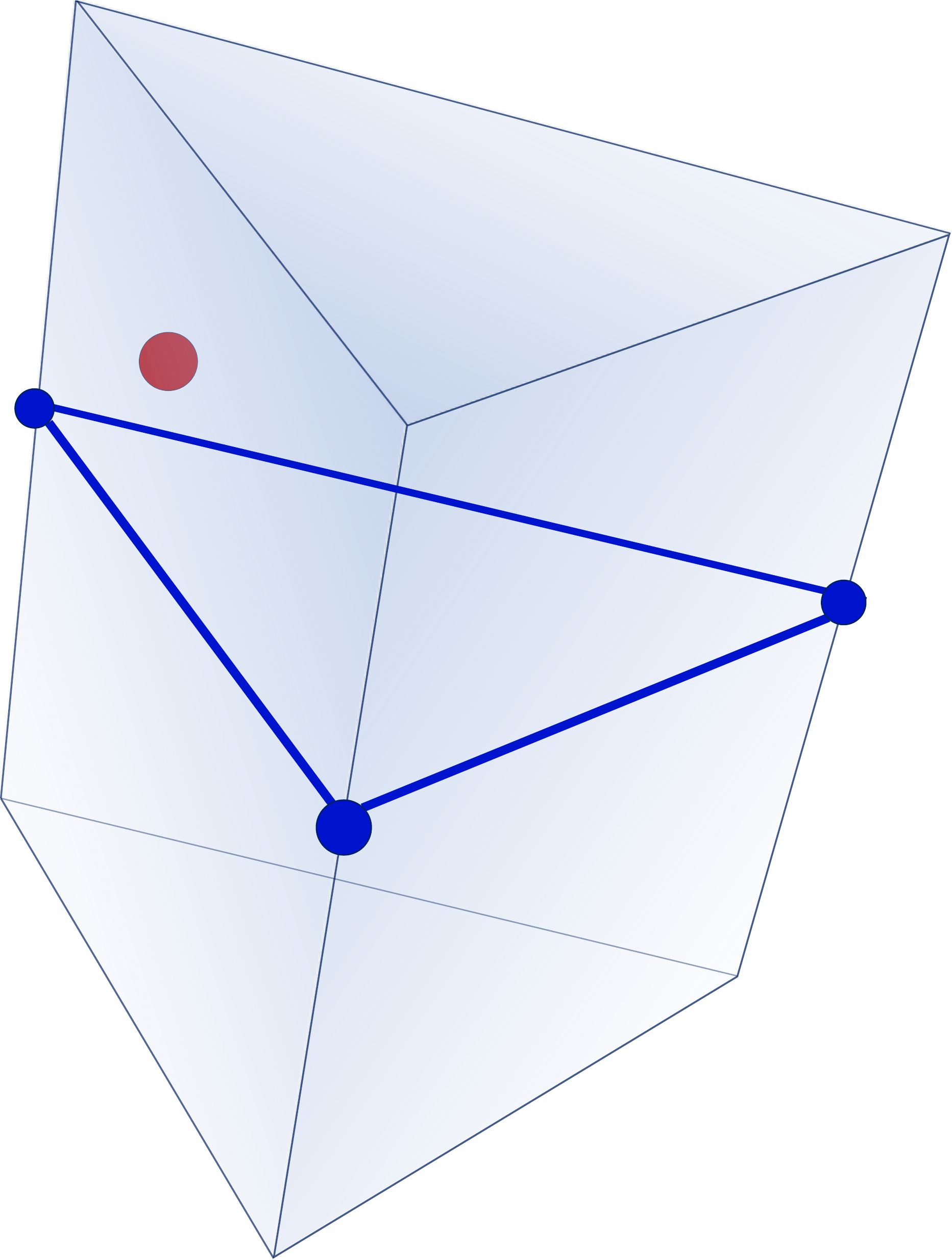}
			\put (30, 120) {Case $1$}
			\put (35, 49) {$e_j$}
			\put (9, 76) {$\widehat{\vec{X}}_i$}
			\end{overpic}
			\end{subfigure}
			\hspace{-1cm}
			\begin{subfigure}[b]{0.3\textwidth}
			\begin{overpic}[width=1.1\textwidth]{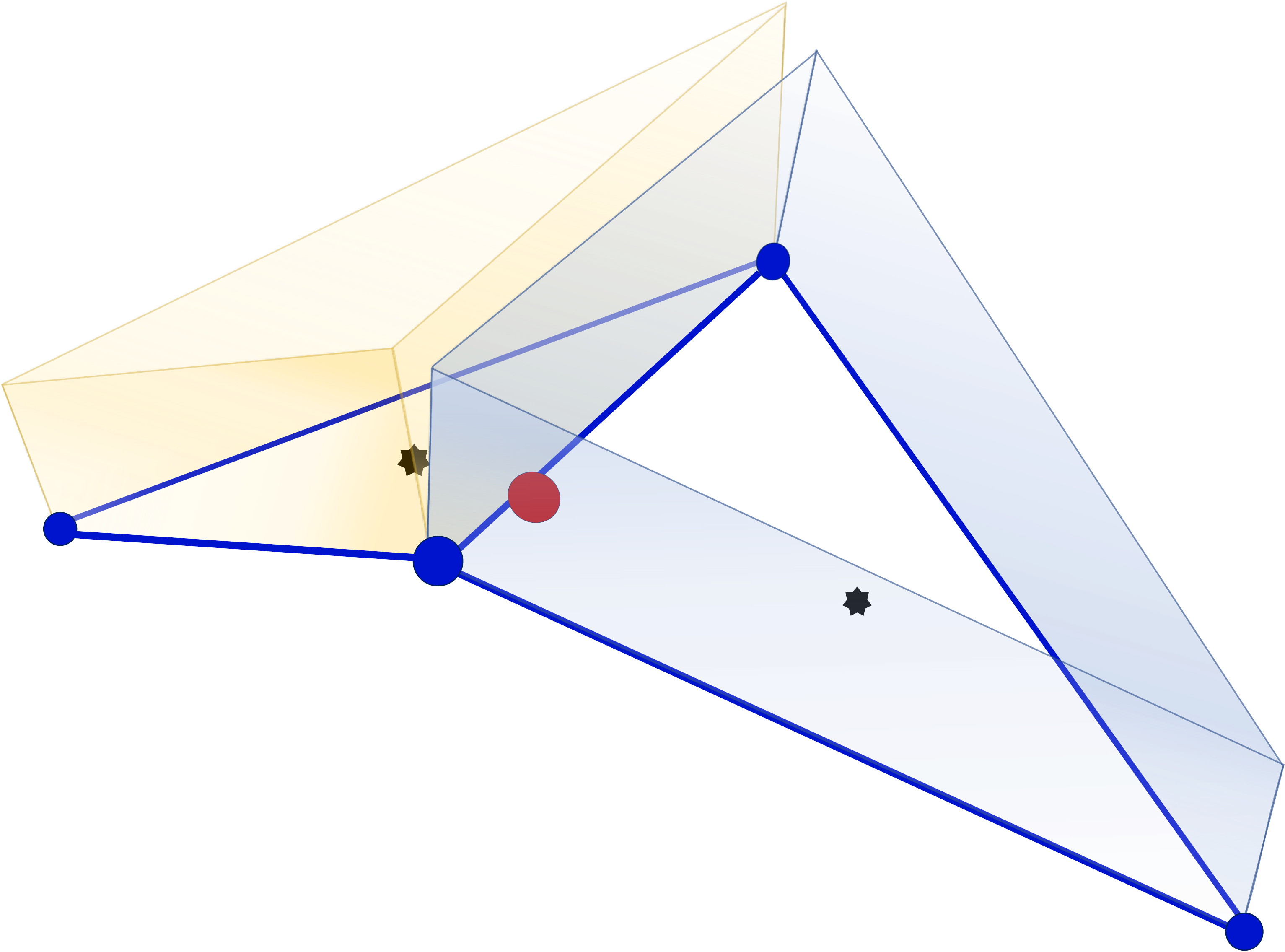}
			\put (40, 101) {Case $2$}
			\put (18, 35) {$e_{k_1}$}
			\put (45, 31.5) {$\widehat{\vec{X}}_i$}
			\put (69, 20) {$e_{k_2}$}
			\put (53.5, 40.5) {$k$}
			\end{overpic}
			\end{subfigure}
			\hspace{0.8cm}
			\begin{subfigure}[b]{0.3\textwidth}
			\begin{overpic}[width=1.1\textwidth]{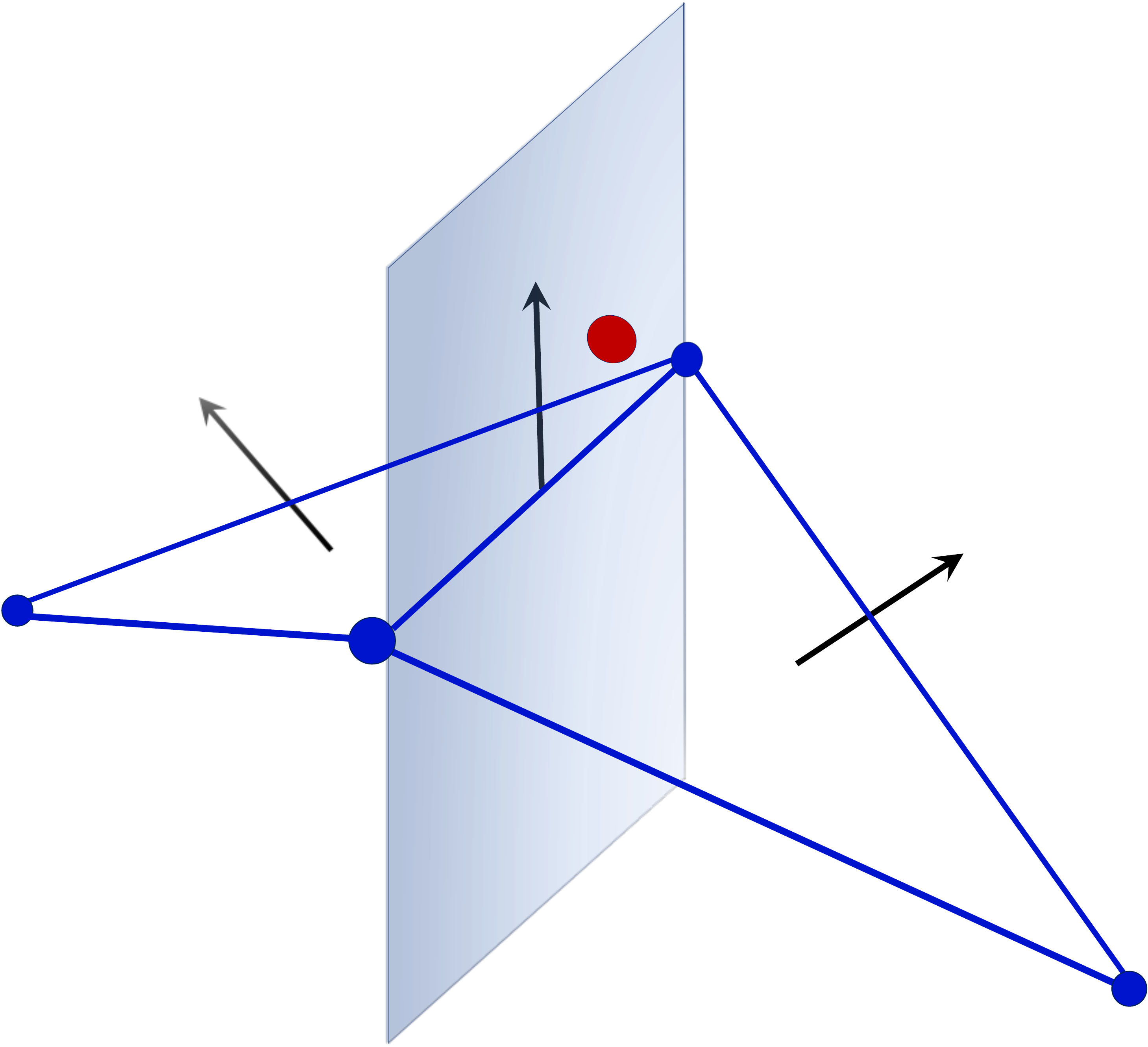}
			\put (35, 101) {Case $3$}
			\put (53, 78) {$\mathbf{P}$}
			\put (18, 40) {$e_{k_1}$}
			\put (69, 28) {$e_{k_2}$}
			\put (55, 65) {$\widehat{\vec{X}}_i$}
			\put (75, 47) {$\vec{n}_{k_1}$}
			\put (5, 55) {$\vec{n}_{k_2}$}
			\put (52, 44) {$k$}
			\put (39, 63) {$\overline{\vec{n}}$}
			\end{overpic}
			\end{subfigure}
			\vspace{0.5cm}
			\caption{\setstretch{1.4}\textbf{Examples for each case of the nearest element search (Algorithm \ref{alg:FindingNearestElement}).} 
			\textbf{Case 1:} Node $i$ (located at $\widehat{\vec{X}}_i$) is contained within the triangular prism defined by the element with the closest centroid ($e_j$). \textbf{Case 2:} Node $i$ resides in either prism defined by the elements (Elements $e_{k_1}$ and $e_{k_2}$) associated with the nearest edge, $k$, but not necessarily the element with the nearest centroid ($*$). \textbf{Case 3:} Node $i$ exists in the region between the two prisms, thus a plane, $\mathbf{P}$, is constructed to determine on which side of Edge $k$ Node $i$ exists. Here $\vec{n}_{k_1}$ and $\vec{n}_{k_2}$ define the unit normals for Elements $e_{k_1}$ and $e_{k_2}$, respectively, and $\overline{\vec{n}} =0.5[\vec{n}_{k_1}+\vec{n}_{k_2}]$. }\label{fig:NearestElementSearch} 
\end{figure}
The new mesh, $\widehat{\mathbf{M}}$, is  an improved approximation of the 3D surface defined by $\Omega_t$, and therefore by definition will not co-exist with $\mathbf{M}$. As such, we can not simply assume that all of the nodes in $\widehat{\mathbf{N}}$ will exist within elements of $\mathbf{M}$. Here we describe a method to find the nearest element in $\mathbf{M}$ for each node in $\widehat{\mathbf{N}}$ (Algorithm \ref{alg:FindingNearestElement}). This process is divided into three successive cases. \textbf{Case 1:} For each node, $i$, we begin by locating the nearest element centroid in $\mathbf{M}$ (Element $e_j$). Element $e_j$ is accepted to be the nearest element if Node $i$ exists within the triangular prism defined by the element (Figure \ref{fig:NearestElementSearch}, Case 1). \textbf{Case 2:} If the test in Case 1 fails then we locate more potential nearest elements by identifying the nearest edge, Edge $k$, in $\mathbf{M}$ (Figure \ref{fig:NearestElementSearch}, Case 2). Edge $k$ is associated with two elements,  $e_{k_1}$ and $e_{k_2}$ (external edges have only one associated element so we choose this element as the nearest element).
If  Node $i$ is contained, in the prism of either $e_{k_1}$ or $e_{k_2}$, then we accept the relevant element. In Figure \ref{fig:NearestElementSearch}, Case 2, it is clear that the centroid (black star) of Element $e_{k_1}$ is the nearest centroid to Node $i$, however Node $i$ resides in the prism defined by Element $e_{k_2}$. 
If both Element $e_{k_1}$ and Element $e_{k_2}$ are rejected, then we need to consider a third case. \textbf{Case 3:} the node lies in the region between the two prisms defined by the elements.
To establish whether Element $e_{k_1}$ or Element $e_{k_2}$ is nearest to Node $i$, a plane, $\mathbf{P}$, dividing this region is constructed (Figure \ref{fig:NearestElementSearch}, Case 3), and we determine on which side of the plane Node $i$ exists, thus providing the nearest element. Plane $\mathbf{P}$ is defined as the plane containing both Edge $k$ and the mean of the associated element normals, $\overline{\vec{n}} =0.5[\vec{n}_{k_1}+\vec{n}_{k_2}]$.  Here $\vec{n}_{k_1}$ and $\vec{n}_{k_2}$ are the unit normals of elements $e_{k_1}$ and $e_{k_2}$, respectively. The situation may arise where Node $i$ is contained in both prisms defined by elements $e_{k_1}$ and $e_{k_2}$. In this case one should follow the process outlined for Case 3. The steps of this process are outlined in Algorithm 1, and the cases are illustrated in Figure \ref{fig:NearestElementSearch}. Over a mesh with sufficient curvature, many element prisms will overlap, and Node $i$ can be contained in multiple prisms. As such, it is important to ensure the selected element prism is also the nearest, as determined using the centroids and edges. The procedure to find the nearest element collapses down for 2D domains, where one simply needs to find the triangular element containing Node $i$ without concerns about the third dimension.
\begin{algorithm}[t]
		\setstretch{0.95}
		\SetAlgoLined
		 \For{Each Node $i$ in $\widehat{\mathbf{M}}$}{
				Find the nearest element centroid in $\mathbf{M}$, Element $j$\;
			\eIf{Node $i$ is contained within the prism defined by Element $j$}{
			 Return Element $j$\; }{
			 Find the nearest edge in $\mathbf{M}$, Edge $k$\;
					 \uIf{Node $i$ is contained within the prism defined by either element associated to Edge $k$, Element $e_{k_1}$ or $e_{k_2}$}{
					 Return relevant element, $e_{k_1}$ or $e_{k_2}$\;  }
					 \ElseIf{Node $i$ is contained within either both, or neither, prisms defined by either elements, $e_{k_1}$ or $e_{k_2}$}{
						Determine if node $i$ is above or below the plane, $\mathbf{P}$, containing Edge $k$ and $\overline{\vec{n}} =0.5[\vec{n}_{k_1}+\vec{n}_{k_2}]$, and orientated towards $e_{k_1}$, where $\vec{n}_{k_1}$ and $\vec{n}_{k_2}$ are the unit normals associated with elements $e_{k_1}$ and $e_{k_2}$, respectively.\;
					 \eIf{Above $\mathbf{P}$}{
					 Return Element $e_{k_1}$\;  }{
					 Return Element $e_{k_2}$\; }}}}
		 \caption{Finding the nearest element in $\mathbf{M}$ for each node in $\widehat{\mathbf{M}}$.}\label{alg:FindingNearestElement}
		\end{algorithm}

\section{Results}
	 In this section we show that the history-dependent remeshing technique described in this paper provides consistent results for spatial deformation and strain energy across the domain after both 2D and 3D deformation. To do so, we examine the spatial error, $E_{\vec{x}}^i$, and strain error, $E_{\Psi}^j$, across $\widehat{\mathbf{M}}$ that arises following remeshing. Additionally, we consider how this history-dependent remeshing technique improves simulation stability across the capillary plexus deformation shown in Figure \ref{fig:ExampleOfMeshDegredation}.

	\subsection{Spatial information is preserved after history-dependent remeshing}
		To examine the spatial error across $\widehat{\mathbf{M}}$ following remeshing, we consider a set of 2D meshes discretising $\Omega_0$ ($3\times 3$ square, Figure \ref{fig:OldAndNewMesh} $\mathbf{M}_0$) with increasing mesh refinement of $\mathbf{M}_0$. Here we define refinement using the initial target edge length at mesh generation. Each mesh is deformed, remeshed once using a set of meshes with increasing refinement, and mapped back to the initial geometry, $\Omega_0$, as described above in the computational methodology and shown in Figure \ref{fig:OldAndNewMesh}. The continuous deformation across the mesh is defined at the nodes, and is a $P_1$ variable.  In this example, $\Omega_0$ is deformed using $\vec{f}(\vec{x}_i)=[x_i^2,y_i^2]$, as shown in Figure \ref{fig:OldAndNewMesh}, where $\vec{x}_i$ is the position of Node $i$, however in practice the deformation is not usually known. For each refinement, $\widehat{\pmb{\EuScript{M}}}_0$ is deformed through $\vec{f}(\vec{x})$  and the spatial error 
		 \begin{align*}  
			E_{\vec{x}}^i= ||\widehat{\vec{X}}_i - \vec{f}(\widehat{\vec{x}}_i)||^2,
		\end{align*}
	  is considered between the deformed mapped mesh and the new mesh.\\

	   \begin{figure}[b!]
		\centering
		\begin{subfigure}[b]{0.45\textwidth}
		\hspace{-0.8cm}
		\begin{overpic}[width=\textwidth, trim={1.5cm 11cm 11cm 10cm},clip, left]{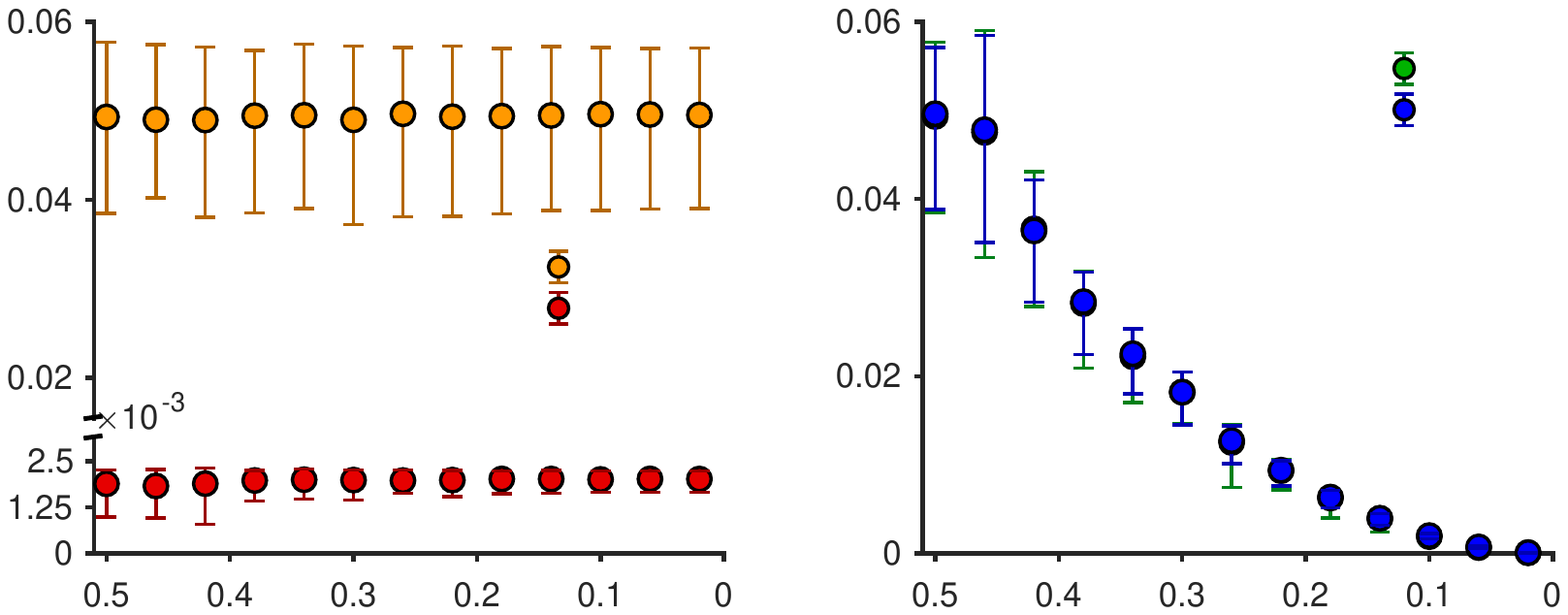}
		\put (1, 90) {$a)$}
		\put (44, -1) {\fontsize{9}{9}\selectfont Edge length}
		\put(4, 40){\rotatebox{90}{ \fontsize{9}{9}\selectfont Error}}
		\put (18, 90) {\fontsize{10}{10}\selectfont Spatial error with increasing remeshed}
		\put (25, 84) {\fontsize{10}{10}\selectfont mesh refinement over 2D mesh}
		\put (78.5,47.1) {\fontsize{8}{8}\selectfont Coarse $\mathbf{M}_0$}
		\put (78.5, 42.1) {\fontsize{8}{8}\selectfont Fine $\mathbf{M}_0$}
		\put(80.65, 11.3){\color{black}\circle*{1.8}}
		\put(80.65, 11.3){\color{blue}\circle*{1.4}}
		\put(21.25, 11.3){\color{black}\circle*{1.8}}
		\put(21.25, 11.3){\color{green}\circle*{1.4}}
		\end{overpic}
		\end{subfigure}
		\hspace{0.2cm}
		\begin{subfigure}[b]{0.45\textwidth}
		\begin{overpic}[width=\textwidth, trim={10.5cm 11cm 2cm 10cm},clip, right]{Fig5.pdf}
		\put (3, 90) {$b)$}
		\put(0, 40){\rotatebox{90}{ \fontsize{9}{9}\selectfont Error}}
		\put (77.7, 71) {\fontsize{8}{8}\selectfont Coarse $\widehat{\mathbf{M}}$}
		\put (77.7, 66.0) {\fontsize{8}{8}\selectfont Fine $\widehat{\mathbf{M}}$}
		\put (44, -1) {\fontsize{9}{9}\selectfont Edge length}
		\put (20, 90) {\fontsize{10}{10}\selectfont Spatial error with increasing initial}
		\put (25, 84) {\fontsize{10}{10}\selectfont mesh refinement over 2D mesh}
		\put(77.05, 11.3){\color{black}\circle*{1.8}}
		\put(77.05, 11.3){\color{red}\circle*{1.4}}
		\put(17.55, 11.3){\color{black}\circle*{1.8}}
		\put(17.55, 11.3){\color{orange}\circle*{1.4}}
		\end{overpic}
		\end{subfigure}
		\vspace{0.1cm}
			\caption{ \setstretch{1.4} \textbf{The median ($\pm$ IQR) spatial error across the deformed square shown in  Figure \ref{fig:OldAndNewMesh} following history-dependent remeshing}. $a)$  The spatial error is unchanged with increasing refinement of the $\widehat{\mathbf{M}}$ remeshed from either a coarse initial mesh (yellow) or a fine initial mesh (red). $b)$ The spatial error decreases with increasing refinement of the initial mesh with both a coarse $\widehat{\mathbf{M}}$ (green) and a fine $\widehat{\mathbf{M}}$ (blue), with only a slight difference between each. The green and blue dots below the $x$-axis in $a)$ denoted the refinement of the initial mesh, either coarse ($0.5$, green) or fine ($0.1$, blue). Similarly the orange and red dots below the $x$-axis in $b)$ denote the refinement of the remeshed mesh, either coarse ($0.5$, orange) or fine ($0.1$, red).\vspace{0.5cm}
			}\label{fig:2DErrorAnal} 
	\end{figure}

Figure \ref{fig:2DErrorAnal} shows the median spatial error  ($\pm$ interquartile range, IQR) over $\widehat{\mathbf{M}}$ for increasing refinement of either $\widehat{\mathbf{M}}$ (a) or $\mathbf{M}_0$ (b). In Figure \ref{fig:2DErrorAnal}a, the refinement of $\widehat{\mathbf{M}}$ increases, while the refinement of $\mathbf{M}_0$ held constant, being either coarse (yellow), or fine (red). 
The spatial error decreases by approximately two orders of magnitude when remeshing from a fine initial mesh, compared to remeshing from a coarse initial mesh. This is confirmed in Figure \ref{fig:2DErrorAnal}b where the refinement of the adapted mesh is held constant, either coarse (green), or fine (blue), while the refinement of the initial mesh is increased. The spatial error decreases with an accompanying narrowing of the IQR as the initial mesh refinement increases. There is little to no additional accuracy attained using a fine  $\widehat{\mathbf{M}}$ over a course one. In practice, this means that the refinement of the new mesh has little impact on the on the remeshing error, but increasing the refinement of the original mesh will reduce remeshing error and provide increased model accuracy. Note that in these (and all subsequent) plots the minimum and maximum error follow the same pattern as the IQR. \\

 \begin{figure}[t!]
	\begin{subfigure}[t]{0.45\textwidth}
 \hspace{0.4cm}
		\begin{overpic}[width=0.8\textwidth, trim={-9cm 0cm 60cm 10cm},clip, left]{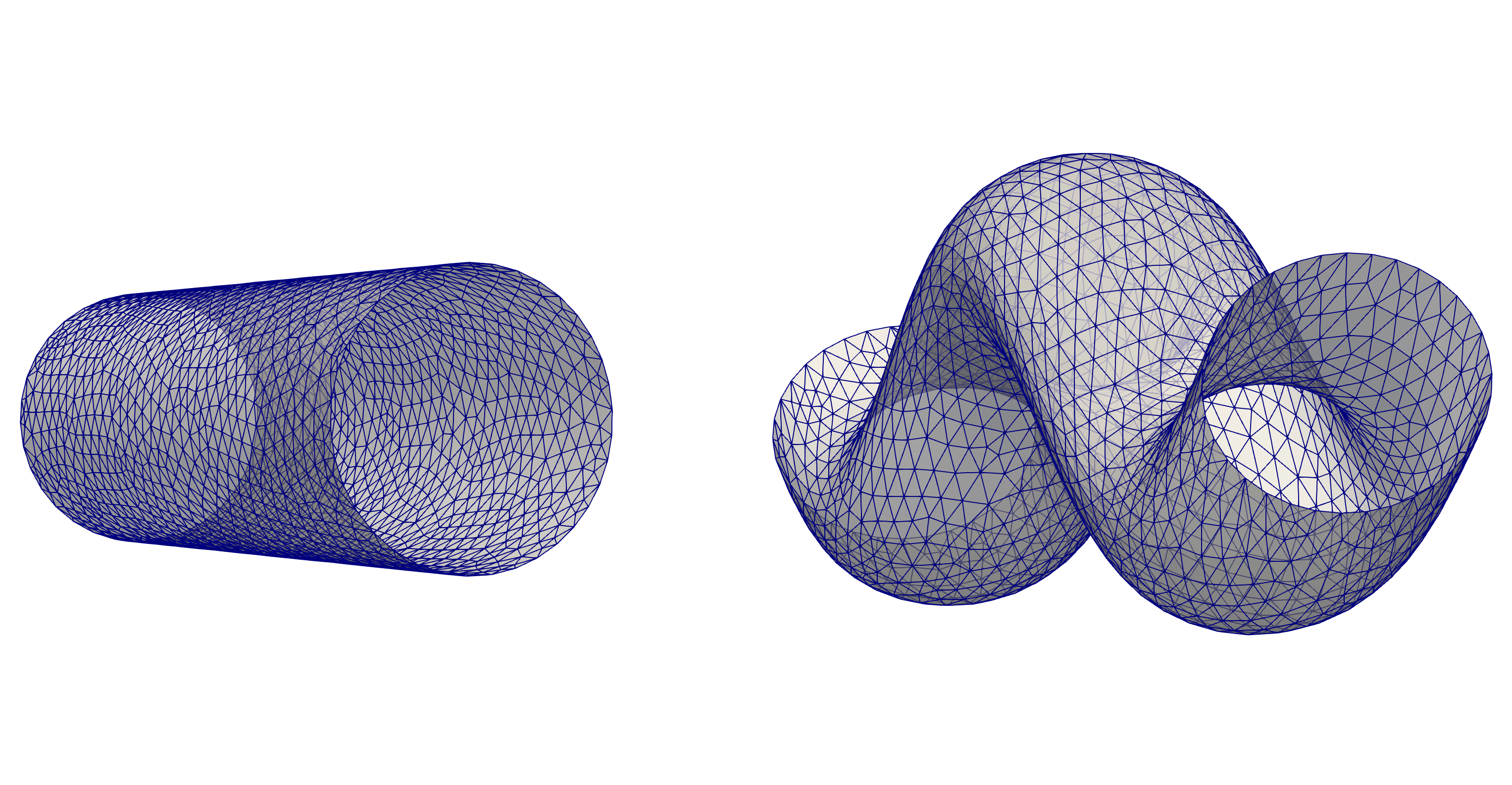}
		\put (-23, 63) {$a)$}
		\put (6 , 15) {$t_0$}
		\put (90, 35) {$\mathrel{\contour{black}{${\longrightarrow}$}}$}
				\put (90, 35) {$\mathrel{\contour{black}{${\longrightarrow}$}}$}
		\put (90, 40) {\fontsize{9}{9}\selectfont $\vec{g}(\vec{x})$}
		\end{overpic}
	\end{subfigure}
\hspace{-0.5cm}
	\begin{subfigure}[t]{0.45\textwidth}
		\begin{overpic}[width=0.8\textwidth, trim={50cm 0cm 1cm 10cm},clip, right]{Fig6A.png}
		\put (13 , 15) {$t_1$}
		\end{overpic}
	\end{subfigure}
	\begin{subfigure}[t]{0.49\textwidth}
		\hspace{-0.55cm}
		\begin{overpic}[width=\textwidth, trim={1.5cm 11cm 11cm 10cm},clip, left]{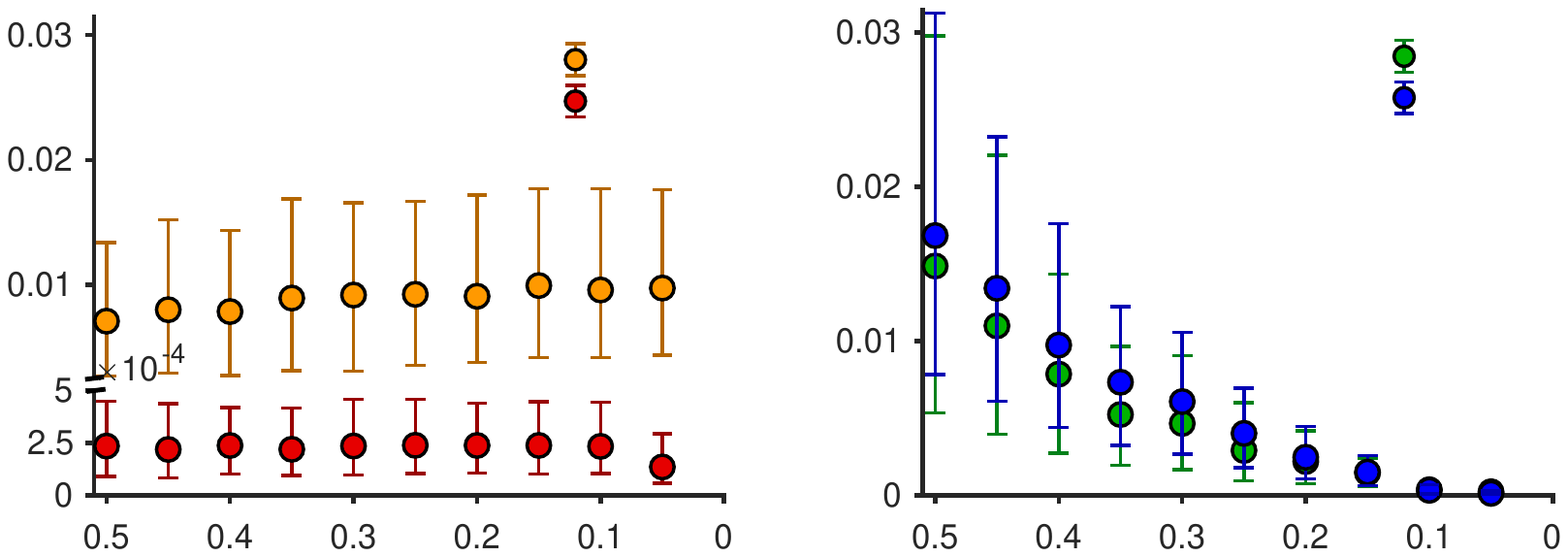}
			\put (3, 87) {$b)$}
			\put (50, 4) {\fontsize{9}{9}\selectfont Edge length}
			\put(2, 41){\rotatebox{90}{ \fontsize{9}{9}\selectfont Error}}
			\put (81, 67.5) {\fontsize{8}{8}\selectfont Coarse $\mathbf{M}_0$}
			\put (81, 62.5) {\fontsize{8}{8}\selectfont Fine $\mathbf{M}_0$}
			\put (19, 86) {\fontsize{10}{10}\selectfont Spatial error with increasing remeshed}
			\put (25, 80) {\fontsize{10}{10}\selectfont mesh refinement over 3D mesh}
			\put(88.15, 14.0){\color{black}\circle*{1.8}}
			\put(88.15, 14.0){\color{blue}\circle*{1.4}}
			\put(36.05, 14.0){\color{black}\circle*{1.8}}
			\put(36.05, 14.0){\color{green}\circle*{1.4}}
		\end{overpic}
	\end{subfigure}
	\begin{subfigure}[t]{0.49\textwidth}
		\begin{overpic}[width=\textwidth, trim={10.5cm 11cm 2cm 10cm},clip, right]{Fig6B.pdf}
			\put (3, 87) {$c)$}
			\put (47, 4) {\fontsize{9}{9}\selectfont Edge length}
			\put(-1.2, 41){\rotatebox{90}{ \fontsize{9}{9}\selectfont Error}}
			\put (77.3, 67.9) {\fontsize{8}{8}\selectfont Coarse $\widehat{\mathbf{M}}$}
			\put (77.3, 63) {\fontsize{8}{8}\selectfont Fine $\widehat{\mathbf{M}}$}
			\put (22, 86) {\fontsize{10}{10}\selectfont Spatial error with increasing initial}
			\put (26, 80) {\fontsize{10}{10}\selectfont mesh refinement over 3D mesh}
			\put(84.5, 14.0){\color{black}\circle*{1.9}}
			\put(84.5, 14.0){\color{red}\circle*{1.5}}
			\put(32.45, 14.0){\color{black}\circle*{1.9}}
			\put(32.45, 14.0){\color{orange}\circle*{1.5}}
		\end{overpic}
	\end{subfigure}
			\caption{\setstretch{1.4} \textbf{The median ($\pm$ IQR) spatial error across a deformed cylinder following history-dependent remeshing.}  $a)$ A cylindrical mesh is deformed by 
			$\vec{g}(\vec{x}_i) = [r_i\cos(\theta_i)+\sin(z_i),~ 1.5r_i\sin(\theta_i),~z_i]$. $b)$  The spatial error across $\widehat{\mathbf{M}}$ is unchanged with increasing mesh refinement of $\widehat{\mathbf{M}}$, remeshed from either a coarse initial mesh (yellow) or a fine initial mesh (red). $c)$ The spatial error across $\widehat{\mathbf{M}}$ (both coarse, green, or fine, blue) reduces sharply with increasing refinement of the initial mesh. The green and blue dots below the $x$-axis in $a)$ denoted the coarse ($0.4$) and fine ($0.05$) refinement of the initial mesh, respectively. Similarly the orange and red dots below the $x$-axis in $b)$ denote the coarse ($0.4$) and fine ($0.05$) refinement of the remeshed mesh, respectively.}
		\label{fig:SpatialErrorInDefomredCylinder}
\end{figure}

 These results are confirmed in 3D, as shown in Figure \ref{fig:SpatialErrorInDefomredCylinder}, where a cylinder is deformed to a larger sinusoidal elliptical cylinder (Figure \ref{fig:SpatialErrorInDefomredCylinder}a), remeshed, and analysed as described above for the 2D example. Here the cylinder is deformed though
 $\vec{g}(\vec{x}_i) = [r_i\cos(\theta_i)+\sin(z_i),~ 1.5r_i\sin(\theta_i),~z_i]$,
 where $\theta_i = \arctan{(y_i/x_i)}$ and $r_i = \sqrt{x_i^2+y_i^2}$. Again, increasing refinement of \nmesh has little effect on the spatial error (Figure \ref{fig:SpatialErrorInDefomredCylinder}b and \ref{fig:SpatialErrorInDefomredCylinder}c). For this 3D example, the spatial error declines with increasing $\mathbf{M}_0$ refinement. Increasing the refinement of $\mathbf{M}_0$ is also accompanied by a narrowing of the standard deviation, as seen in 2D (Figure \ref{fig:2DErrorAnal}b). These results suggest that the spatial accuracy of the history-dependent remeshing technique is dependent only on the refinement of the initial mesh. This is unsurprising as additional spatial information can not be attained mid-way through a simulation by re-discretising the deformed geometry. 
 Rather, the spatial accuracy of the original mesh limits the attainable information about the initial configuration. 
 
\subsection{The hyperelastic strain energy density across $\widehat{\mathbf{M}}$ reproduces that of $\mathbf{M}$}\label{StrainAnalysisSection}
We now examine the error introduced to the strain across the deformed, $\mathbf{M}$, and remeshed, $\widehat{\mathbf{M}}$, meshes. The strain, a $P_0$ variable, is defined on and is continuous over each element, but is discontinuous over the domain. Here we assess the error in the hyperelastic strain energy density subsequent to history-dependent remeshing to understand how simulations behave following this technique. We also consider error minimisation by optimising the initial, and adapted mesh refinement and the remeshing frequency. We model the deforming square domain, shown in Figure \ref{fig:OldAndNewMesh}, as a hyperelastic membrane, where the stress-strain relationship can be derived from a strain energy density function, $\Psi$. Following Kruger \etal \cite{kruger2011efficient} we adopt the Skalak strain energy density function \cite{skalak1973},
\begin{equation} \label{SkakalStrainEnergy}
	\Psi(I_1,I_2) = \frac{\kappa_s}{12}(I_1^2+2I_1-2I_2)+ \frac{\kappa_{\alpha}}{12}I_2^2,
\end{equation}
where $\kappa_{\alpha}$ and $\kappa_s$ are the area and strain dilation moduli respectively, and $I_1$ and $I_2$ are the 2D principal strain invariants, describing the strain, and the dilation state of the membrane, respectively. 
The Skalak strain energy density function was first developed by Skalak \cite{skalak1973} in 1973 to describe the deformation of red blood cells, and has been since used in to model red blood cell movement \cite{kruger2011efficient}. We assume each element is both mechanically isotropic and homogeneous, and that deformation is linear across each element. For further details on the numerical implementation of this model see Kruger 
\etal \cite{kruger2011efficient}. \\

 \begin{figure}[h!]
	\begin{subfigure}[b]{0.49\textwidth}
		\begin{overpic}[width=1\textwidth, trim={1cm 5cm 0cm 7cm},clip]{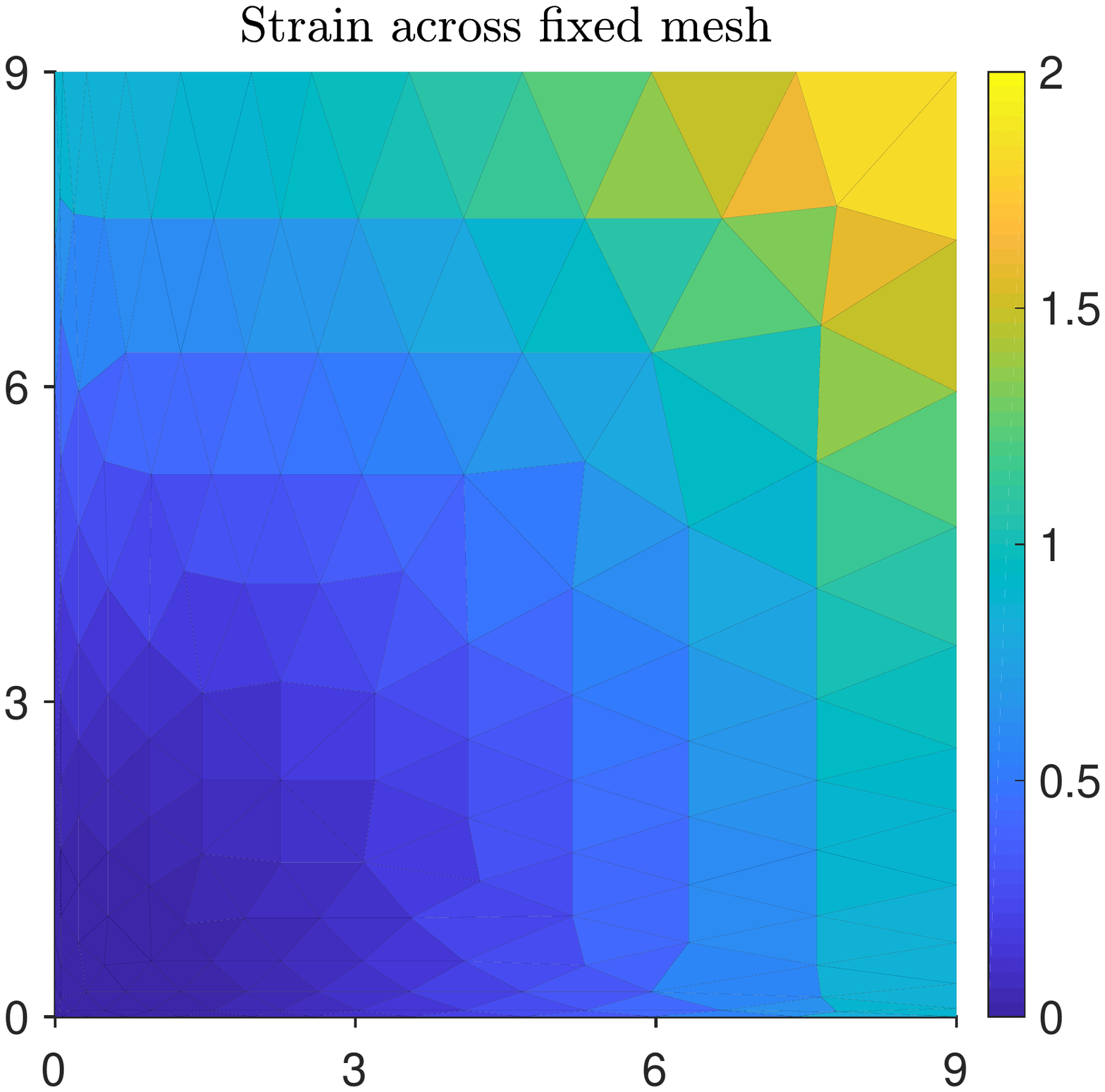}
	\put (-4, 94) {$a)$}
	\put (22, 95) {\fontsize{10}{10}\selectfont Strain across fixed mesh}
		\end{overpic}
		\vspace{-1.2cm}
	\end{subfigure}
	\begin{subfigure}[b]{0.49\textwidth}
		\begin{overpic}[width=1\textwidth, trim={1cm 5cm 0cm 7cm},clip]{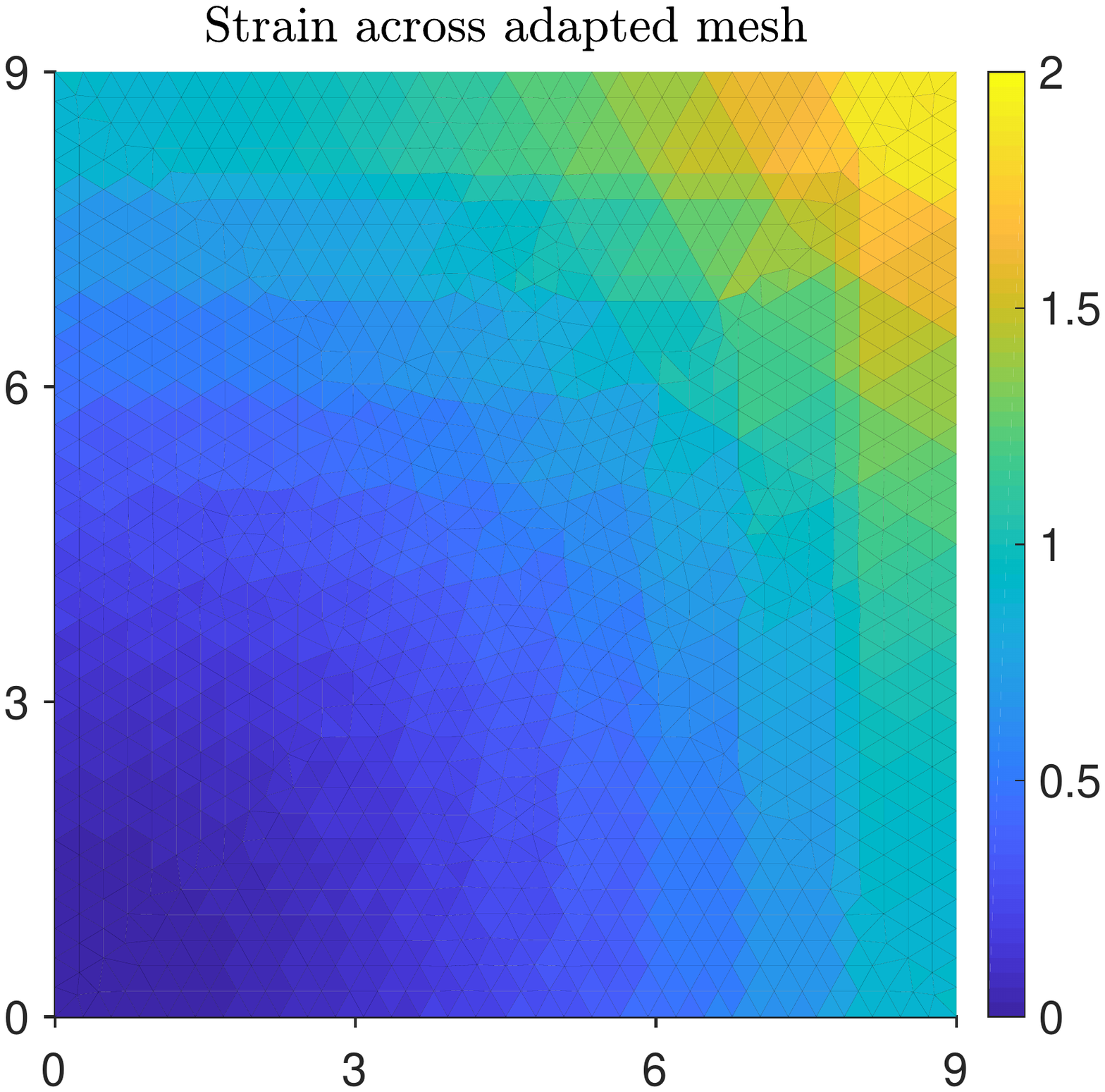}
	\put (-4, 94) {$b)$}
	\put (18, 95) {\fontsize{10}{10}\selectfont Strain across remeshed mesh}
		\end{overpic}
		\vspace{-1.2cm}
	\end{subfigure}
	\begin{subfigure}[b]{0.49\textwidth}
	\vspace{0.2cm}
	\hspace{-0.8cm}
		\begin{overpic}[width=\textwidth, trim={1cm 9cm 10.5cm 10cm},clip,left]{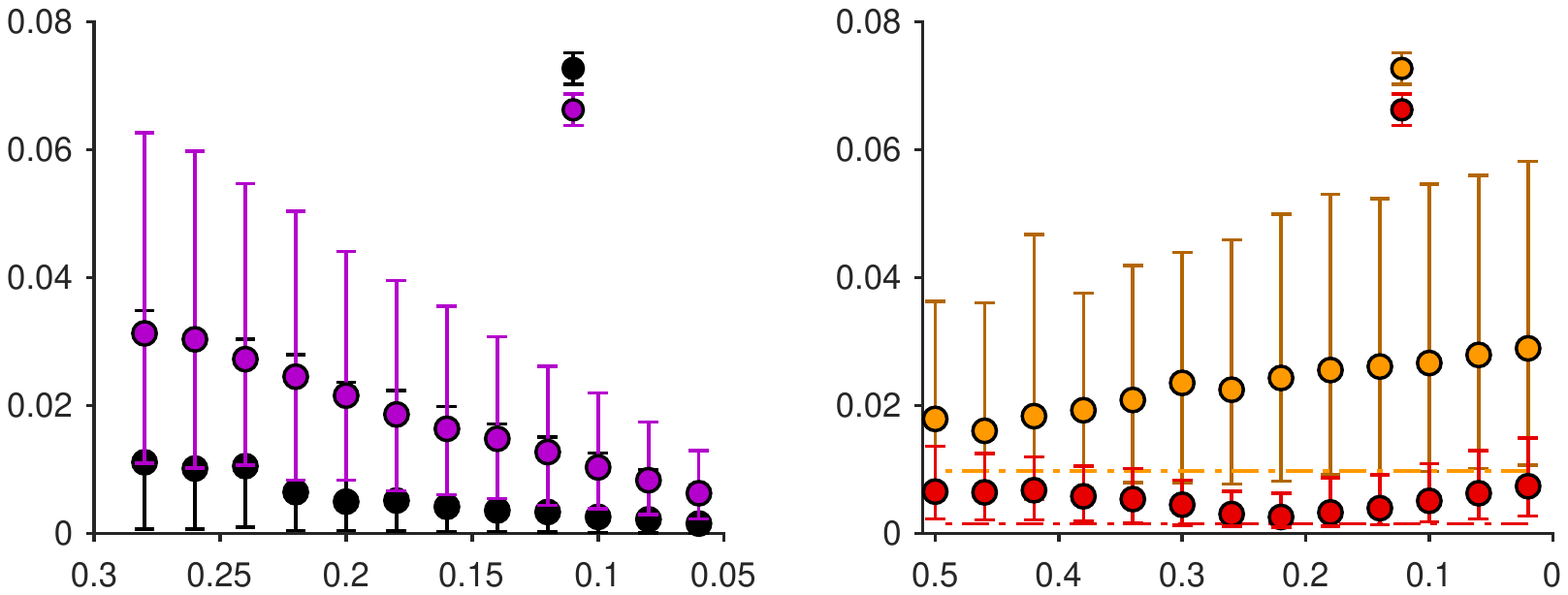}
	\put (5, 88) {$c)$}
			\put (71, 75.5) {\fontsize{8}{8}\selectfont $\mathbf{M}$}
			\put (71, 70.9) {\fontsize{8}{8}\selectfont $\widehat{\mathbf{M}}$}
			\put (42,20) {\fontsize{9}{9}\selectfont Edge length}
			\put(6.4, 47){\rotatebox{90}{ \fontsize{9}{9}\selectfont Error}}
		\put (26, 89) {\fontsize{9}{9}\selectfont Strain error with increasing}
		\put (35, 84) {\fontsize{9}{9}\selectfont mesh refinement}
		\end{overpic}
	\end{subfigure}
	 \vspace{-0.2cm}
	\begin{subfigure}[b]{0.49\textwidth}
		\begin{overpic}[width=\textwidth, trim={10.5cm 9cm 1cm 10cm},clip]{Fig7C.pdf}
	\put (-3, 88) {$d)$}
			\put (63, 75.4) {\fontsize{8}{8}\selectfont Coarse $\mathbf{M}_0$}
			\put (63, 71.2) {\fontsize{8}{8}\selectfont Fine $\mathbf{M}_0$}
			\put (35, 20) {\fontsize{9}{9}\selectfont Edge length}
			\put(-2, 47){\rotatebox{90}{ \fontsize{9}{9}\selectfont Error}}
		\put (22, 89) {\fontsize{9}{9}\selectfont Strain error with increasing}
		\put (23, 84) {\fontsize{9}{9}\selectfont remeshed mesh refinement}
			\put(45, 29){\color{black}\circle*{1.4}}
			\put(45, 29){\color{green}\circle*{1.0}}
			\put(69.5, 29){\color{black}\circle*{1.4}}
			\put(69.5, 29){\color{lblue}\circle*{1.0}}
		\end{overpic}
	\end{subfigure}
	\vspace{-1.3cm}
	\caption{\setstretch{1.2} \textbf{The median ($\pm$ IQR) strain error across $\widehat{\mathbf{M}}$ decreases with increasing initial mesh refinement.} The Skalak strain energy density over the original mesh, $\mathbf{M}$, and the new mesh, $\widehat{\mathbf{M}}$,  following deformation ($\vec{f}(\vec{x}_i) = [x_i^2,y_i^2]$, Figure \ref{fig:OldAndNewMesh}) are shown in a) and b), respectively. c) The median ($\pm$ IQR) strain error across $\mathbf{M}$ (black) and the successive $\widehat{\mathbf{M}}$ (purple) decreases with increasing mesh refinement. When remeshing, $\widehat{\mathbf{M}}$ is ascribed the same edge length as its preceding $\mathbf{M}$, and as such the error of the associated meshes can be seen by comparing directly for each refinement level. 
	d) The median ($\pm$ IQR) strain error across $\widehat{\mathbf{M}}$, is reduced for all levels of $\widehat{\mathbf{M}}$ refinement if the preceding $\mathbf{M}$ is fine (red), 
	rather than coarse (orange).     
For clarity, the preceding $\mathbf{M}$ is identified with either a green circle (coarse, edge length 0.26) or a blue circle (fine, edge length 0.06), and the median strain error over the coarse $\mathbf{M}$ and the fine $\mathbf{M}$ are shown by the orange and red dashed lines, respectively. }\label{fig:StrainError}
\end{figure}

Figure \ref{fig:StrainError} shows the Skalak strain energy density across $\mathbf{M}$ ((a), initial edge length $0.5$), and $\widehat{\mathbf{M}}$ ((b), initial edge length $0.2$). Here the domain is remeshed only once following deformation, and we use $ \kappa_s = 0.01$ and $\kappa_{\alpha} = 0.01\times10^{-4}$. 
The coarseness of the initial mesh is reflected in the Skalak strain energy density, which shows a pixelated, or granular, nature between elements (Figure \ref{fig:StrainError}a). This granularity is inherited by the adapted mesh, where the echo of the original strain energy density can be seen. This echo is particularly evident in the top right corner of the adapted mesh (Figure \ref{fig:StrainError}b) where element distortion is greatest and suggests that the refinement of the original mesh, and the remeshing frequency may have a large impact on the accuracy of the final mesh. \\

To examine the error in the Skalak strain energy density, we consider the difference between the exact, $\Psi$, and the numerical, $\hat{\Psi}$ (Equation (\ref{SkakalStrainEnergy})), strain energy density  for each element ($j$),
\begin{align*} \label{StrainErrorEquation}
E^j_{\Psi}= |\Psi^j - \hat{\Psi}^j|.
\end{align*}
Owing to the simplicity of the prescribed deformation, we can derive the exact Skalak strain energy density at each element centroid from the Greens deformation tensor. For the deformation given by $\vec{f}(\vec{x})$ in Figure \ref{fig:OldAndNewMesh}, the strain invariants in Equation (\ref{SkakalStrainEnergy}) can be expressed as,
\begin{align*} 
I_1 &= \sum_{i=1}^2(2\Delta\vec{x}_i^j-1)^2-2,\quad \text{and,}\quad I_2 = \prod_{i=1}^2(2\Delta\vec{x}_i^j-1)^2-1,
\end{align*}
where $\Delta\vec{x}^j$ is the displacement of the centroid of Element $j$. For further details of the implementation and the derivation of the Skalak strain energy density, refer to Kruger \etal \cite{kruger2011efficient} and Skalak \cite{skalak1973}, respectively. As we assume the deformation over each element is linear, the strain energy is constant within an element, and thus we can concisely compare the exact and numerical strain energy density at the centre of the element.\\

From Figure \ref{fig:StrainError}c, we see that the median strain error ($\pm$ IQR) across $\mathbf{M}$ (black) and $\widehat{\mathbf{M}}$ (purple) reduces with increasing mesh refinement. Again, mesh refinement is determined by the edge length at the time of mesh generation. 
Both $\mathbf{M}$, and the successive  $\widehat{\mathbf{M}}$ were ascribed the same edge length at mesh generation, and as such the error of the associated meshes can be seen by comparing directly for each refinement level in Figure \ref{fig:StrainError}c. The median strain error over $\widehat{\mathbf{M}}$ is of the same order of magnitude to that of its predecessor at each level of refinement, however is consistently larger than the strain error over the corresponding initial mesh. 
Adjusting the refinement of $\widehat{\mathbf{M}}$ in isolation has little effect on the strain error, as seen in Figure \ref{fig:StrainError}d, where the initial mesh is fixed, either coarse (orange, edge length 0.26) or fine (red, edge length 0.06), with the refinement of $\widehat{\mathbf{M}}$ increasing. For comparison, the strain error of the initial mesh for the coarse and the fine meshes are shown by the orange and the fine dashed lines, respectively. Again there is a notable decrease in the strain error between the coarse and fine initial mesh for each level of mesh refinement over the adapted mesh. The minimum and maximum strain errors follow the same pattern as the IQR.  

\subsection{Increasing remeshing frequency reduces error}
 The deformation examined above, $\vec{f}(\vec{x}_i) = [ x_i^2,y_i^2]$ shown in Figure \ref{fig:OldAndNewMesh}, produced a heterogeneous three-fold dilation of the initial geometry, and caused significant element distortion, with $25\%$ of elements suffering an aspect ratio below $0.3$ (discussed further below in Section 3.4). Remeshing only once after large element distortion (defined here as an element aspect ratio less than $0.6$) is alone insufficient to increase strain accuracy. Though introducing stability to the simulation, this single remeshing event slightly increases the strain error (Figure \ref{fig:StrainError}). As such, we now examine the spatial error and the strain error over a deforming hyperelastic membrane while remeshing multiple times.\\ 
 
As described above in Section \ref{StrainAnalysisSection}, we model the geometry shown in Figure \ref{fig:OldAndNewMesh} as a deforming hyperelastic membrane. We now remesh the domain at defined intervals (between $0-20$ remeshing events) and discretise the deformation, $\vec{f}(\vec{x}_i) = [x^2_i,y^2_i]$, over time, to  produce the equivalent time dependent deformation function 
\begin{equation}
	\vec{F}(\vec{x}_i,t) = \frac{t}{t_{\text{end}}}[x_{0,i}^2,y_{0,i}^2].
\end{equation}
Here $t$ is time in seconds, $t_{\text{end}} = 60~\text{seconds}$ and is the duration of the simulation, and $(x_{0,i}, y_{0,i})$ is the initial position of Node $i$. As this simulation is run for $60$ seconds, the remeshing intervals are equivalent to remeshing frequencies (per second). Based on the results of previous sections, we use an initial edge length of $0.1$ for both the initial and remeshed meshes.\\

\begin{figure}[t!]
 \vspace{0.1cm}
	\begin{subfigure}[b]{0.45\textwidth}
	\hspace{0.25cm}
	\hspace{-0.6cm}
		\begin{overpic}[width=1\textwidth, trim={11cm 11.5cm 2cm 11.1cm},clip]{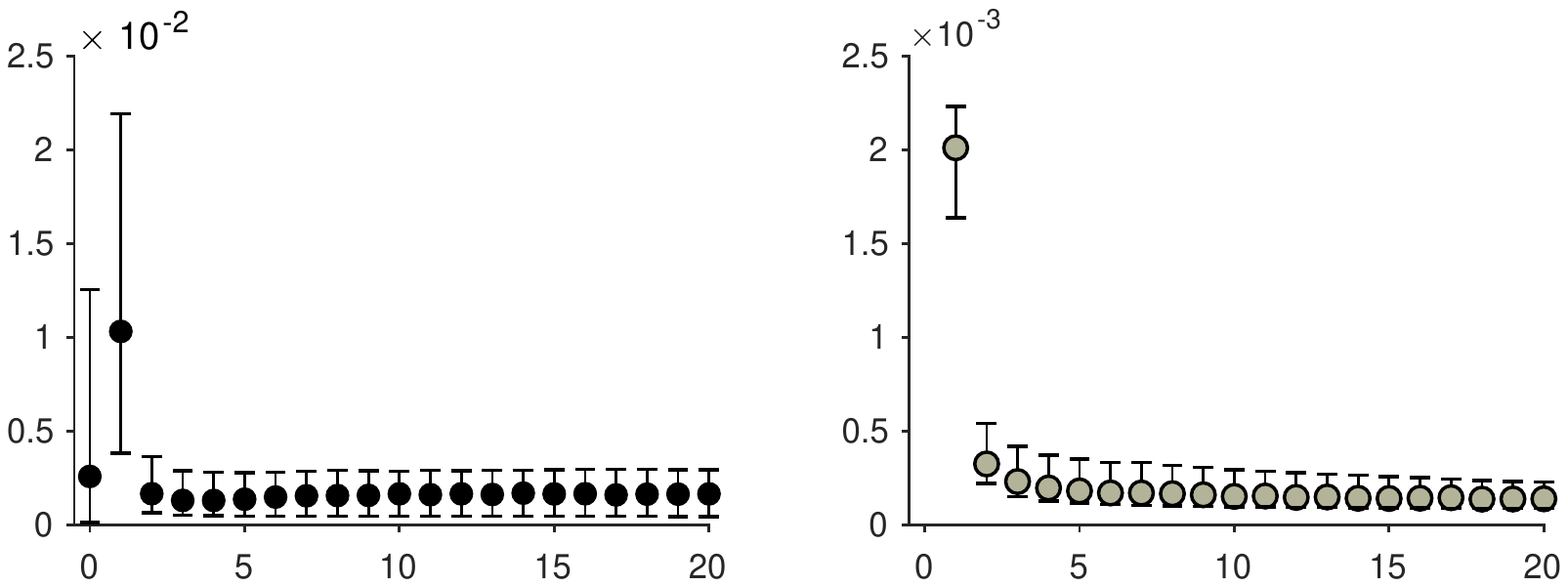}
			\put (35, -1.5) {\fontsize{9}{9}\selectfont Remeshing frequency}
			\put(-4, 27){\rotatebox{90}{ \fontsize{9}{9}\selectfont Spatial error}}
			\put (25, 89) {\fontsize{10}{10}\selectfont Spatial error with increasing} 
			\put (35, 83.5) {\fontsize{10}{10}\selectfont remeshing frequency}
			\put (0, 90) {a)}
		\end{overpic}
	\end{subfigure}
	\begin{subfigure}[b]{0.45\textwidth}
	\hspace{-0.25cm}
		\begin{overpic}[width=1\textwidth, trim={2cm 11.5cm 11cm 11.1cm},clip]{Fig8.pdf}
			\put (39, -1.5) {\fontsize{9}{9}\selectfont Remeshing frequency}
			\put(0, 29){\rotatebox{90}{ \fontsize{9}{9}\selectfont Strain error}}
			\put (29, 89) {\fontsize{10}{10}\selectfont Strain error with increasing} 
			\put (39, 83.5) {\fontsize{10}{10}\selectfont remeshing frequency}
			\put (0, 89) { b)}
		\end{overpic}
	\end{subfigure}
	\vspace{0.4cm}
		\caption{\setstretch{1.4} \textbf{The median ($\pm$ IQR) spatial error (a) and strain error (b) for increasing remeshing frequencies across a square hyperelastic membrane deforming over time.}  
		Here remeshing frequency is measured as remeshing events per minute. }\label{fig:StrainErrorFrequ}
\end{figure}

For each remeshing frequency, we examined the spatial and strain error, as shown in Figure \ref{fig:StrainErrorFrequ}. Note that a frequency of 1 corresponds to only remeshing at the end of the deformation whereas a frequency of $n>1$ corresponds to remeshing at the end and at $n-1$ equally spaced times throughout the deformation. For example a remeshing frequency of 3 corresponds to remeshing (with an initial edge length of $0.1$) at $t=20,~40,$ and $60$ seconds.
The median ($\pm$ IQR) spatial error peaks at a single remeshing event following excessive deformation (Figure \ref{fig:StrainErrorFrequ}a), reflecting the results shown in Figure \ref{fig:SpatialErrorInDefomredCylinder}. Following two remeshing events, the median spatial error falls by an order of magnitude and there is little further accuracy gained for larger remeshing frequencies. Figure \ref{fig:StrainErrorFrequ}b also shows a magnitude drop in the median ($\pm$ IQR) strain error  when the number of remeshing events increases past two. The strain error over the fixed mesh is shown at zero remeshing events, with a larger error and a wider IQR than all other number of remeshing events, with the exception of a single remeshing event. The strain error is increased following a single remeshing event,  reflecting the results in Figure \ref{fig:StrainError}, and again occurring as a consequence of excessive deformation prior to remeshing. For both the spatial error, and the strain error the minimum and maximum error follows the same pattern as the IQR.  \\

 It is necessary to note that the required frequency of remeshing should be guided by the rate and spatial distribution of element distortion. Simulations with more rapid, or with slower element distortion than that shown in the the example here will require more or less frequent remeshing, respectively. We have examined remeshing frequency at constant intervals here for the sake of simplicity, and have shown that after some threshold excessive remeshing adds little further accuracy to the simulation, and that remeshing must occur prior to excessive element deformation.

\subsection{Element aspect ratios over the deformed mesh improves with remeshing at a small expense to the quality of the initial mesh}
Mesh stability is dependent on the quality of the elements forming the mesh. The quality of a given element is often quantified using the element aspect ratio given by $\AR = 2 r_i/r_o$. Here $r_i$ and $r_o$ describe the radii of the inscribed and circumscribed circles, respectively \cite{khan2020surface}. Note that there are several different definitions of the aspect ratio of a triangle, though all have a similar meaning. A higher aspect ratio is indicative of a high quality triangle, while a low aspect ratio is typical of a poor quality triangle for the purpose of a finite element model. An equilateral triangle has an aspect ratio of one, and as a triangle distorts, the aspect ratio falls, approaching zero for a triangle with no surface area. A mesh composed of elements with larger aspect ratios ($\AR \sim 1$) will be more stable under deformation than a mesh deforming with distorted elements (elements with a low aspect ratio).\\

\begin{figure}[t!]
		\centering
		 \vspace{0.5cm}
		 \begin{overpic}[width=0.8\textwidth, trim={2.5cm 8.1cm 1.2cm 8.5cm},clip]{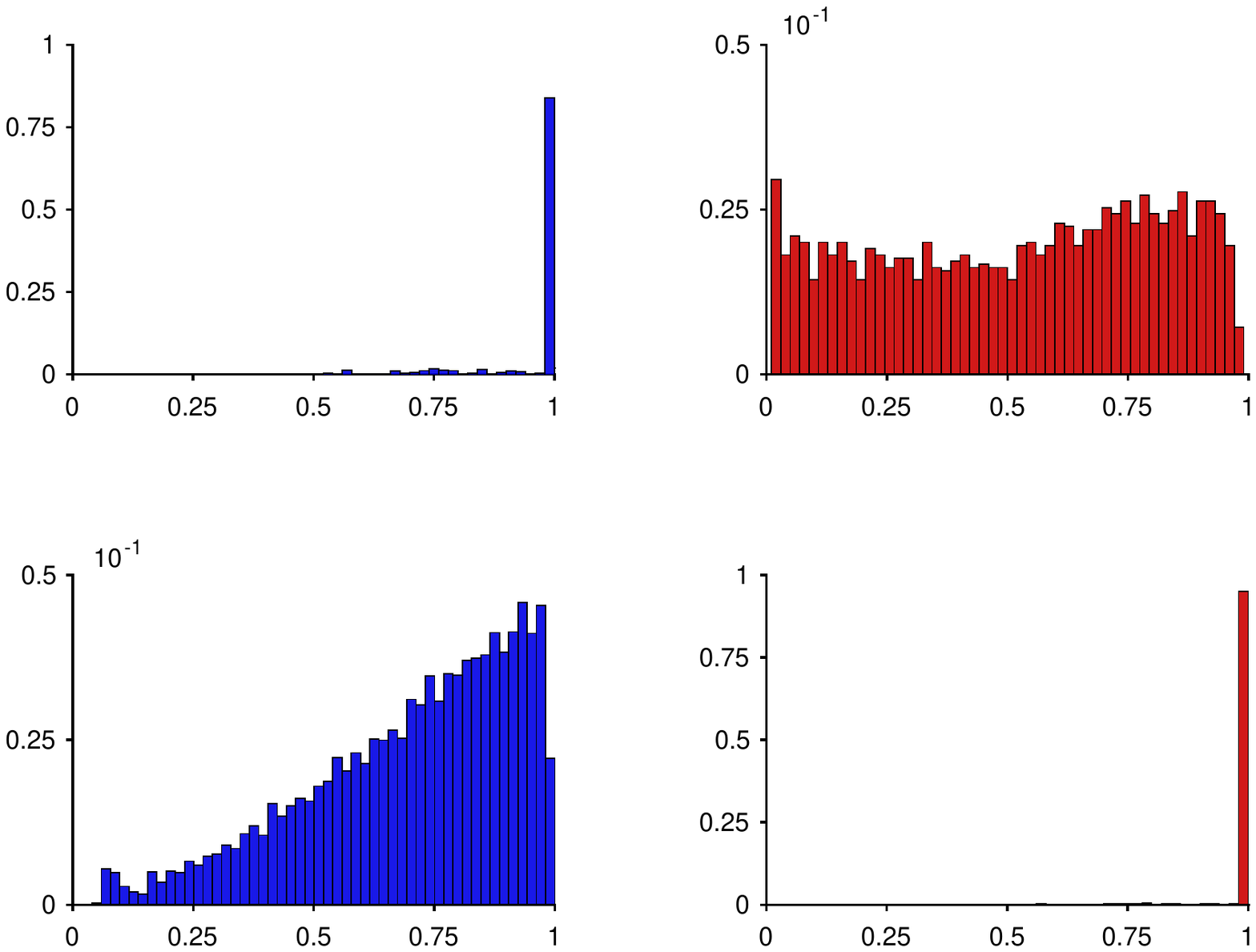}
			 \put (0, 76) {\fontsize{10}{10}\selectfont $a)$}
			 \put (54, 76) {\fontsize{10}{10}\selectfont $b)$}
			 \put (54, 34) {\fontsize{10}{10}\selectfont $d)$}
			 \put (0, 34) {\fontsize{10}{10}\selectfont $c)$}
				\put (21, 74) {\fontsize{10}{10}\selectfont $\mathbf{M}_0$} 
				 \put (76, 74) {\fontsize{10}{10}\selectfont $\mathbf{M}$}
			\put (76, 31) {\fontsize{10}{10}\selectfont$\widehat{\mathbf{M}}$}
				\put (21, 31) {\fontsize{10}{10}\selectfont $\widehat{\pmb{\EuScript{M}}}_0$}
				\put (18, -2.5) {\fontsize{9}{9}\selectfont Aspect ratio}
				\put (73,-2.5) {\fontsize{9}{9}\selectfont Aspect ratio}
				\put (18, 39.5) {\fontsize{9}{9}\selectfont Aspect ratio}
				\put (73, 39.5) {\fontsize{9}{9}\selectfont Aspect ratio}
				\put(-3, 51.2){\rotatebox{90}{ \fontsize{9}{9}\selectfont Distribution}}
				\put(-3, 9){\rotatebox{90}{ \fontsize{9}{9}\selectfont Distribution}}
				\put(52.5, 51.2){\rotatebox{90}{ \fontsize{9}{9}\selectfont Distribution}}
				\put(52.5, 9){\rotatebox{90}{ \fontsize{9}{9}\selectfont Distribution}}
				\end{overpic}
				\vspace{0.3cm}
		 \caption{\setstretch{1.4} \textbf{The element aspect ratio over each mesh.} $\mathbf{M}_0$ (top left), $\mathbf{M}$ (top right), $\widehat{\mathbf{M}}$ (bottom right) and $\widehat{\pmb{\EuScript{M}}}_0$ (bottom left).}. \label{fig:AspectRatios}
	\end{figure}

Remeshing produces a stable mesh with element aspect ratio closer to 1, as shown in Figure \ref{fig:AspectRatios}, thereby increasing the stability of the simulation. Figure \ref{fig:AspectRatios} shows histograms depicting the distribution of the element aspect ratio across the initial undeformed mesh, $\mathbf{M}_0$ ($a$), the deformed mesh, $\mathbf{M}$ ($b$), the remeshed mesh, $\widehat{\mathbf{M}}$ ($d$), and the remeshed mesh mapped to the initial configuration $\widehat{\pmb{\EuScript{M}}}_0$ ($c$). As the original mesh deforms, the elements distort and the median aspect ratio falls from $\AR =1 $ to $\AR=0.56$, and the distribution is relatively uniform, indicating $\mathbf{M}$ is a poor quality mesh (Figure \ref{fig:AspectRatios}b). To recover valid results from the simulation, remeshing is necessary to reinstate a high quality mesh, with a median element aspect ratio close to $\AR =1 $ (Figure \ref{fig:AspectRatios}d). 
This high quality mesh discretising the deformed domain is introduced at a cost to the initial configuration, $\widehat{\pmb{\EuScript{M}}}_0$ (Figure \ref{fig:AspectRatios}c). The aspect ratio across $\widehat{\pmb{\EuScript{M}}}_0$ has a much larger distribution than that of $\widehat{\mathbf{M}}$, with a left skew and a median aspect ratio of  $\AR=0.74$. This left skewed distribution occurs due to the dependence of the mapping on the already distorted $\mathbf{M}$. Importantly, there are less distorted elements across  $\widehat{\pmb{\EuScript{M}}}_0$ compared to $\mathbf{M}$, and the distortion within these elements is less pronounced (Figure \ref{fig:AspectRatios}b and Figure \ref{fig:AspectRatios}c). The presence of highly distorted elements over $\widehat{\pmb{\EuScript{M}}}_0$ contributes to the spatial error and strain error. Increasing remeshing frequencies, or remeshing prior to excessive distortion will reduce the variation in the element aspect ratio over $\widehat{\pmb{\EuScript{M}}}_0$, thereby reducing the spatial error and strain error, as seen in Figure \ref{fig:StrainErrorFrequ}. 

\begin{figure}[t!]
\begin{subfigure}[b]{0.56\textwidth}
 \vspace{4cm}
 \hspace{-0.2cm}
\begin{overpic}[trim = {2cm 3cm 50cm 3cm} , width=0.55\textwidth,  clip=true, left]{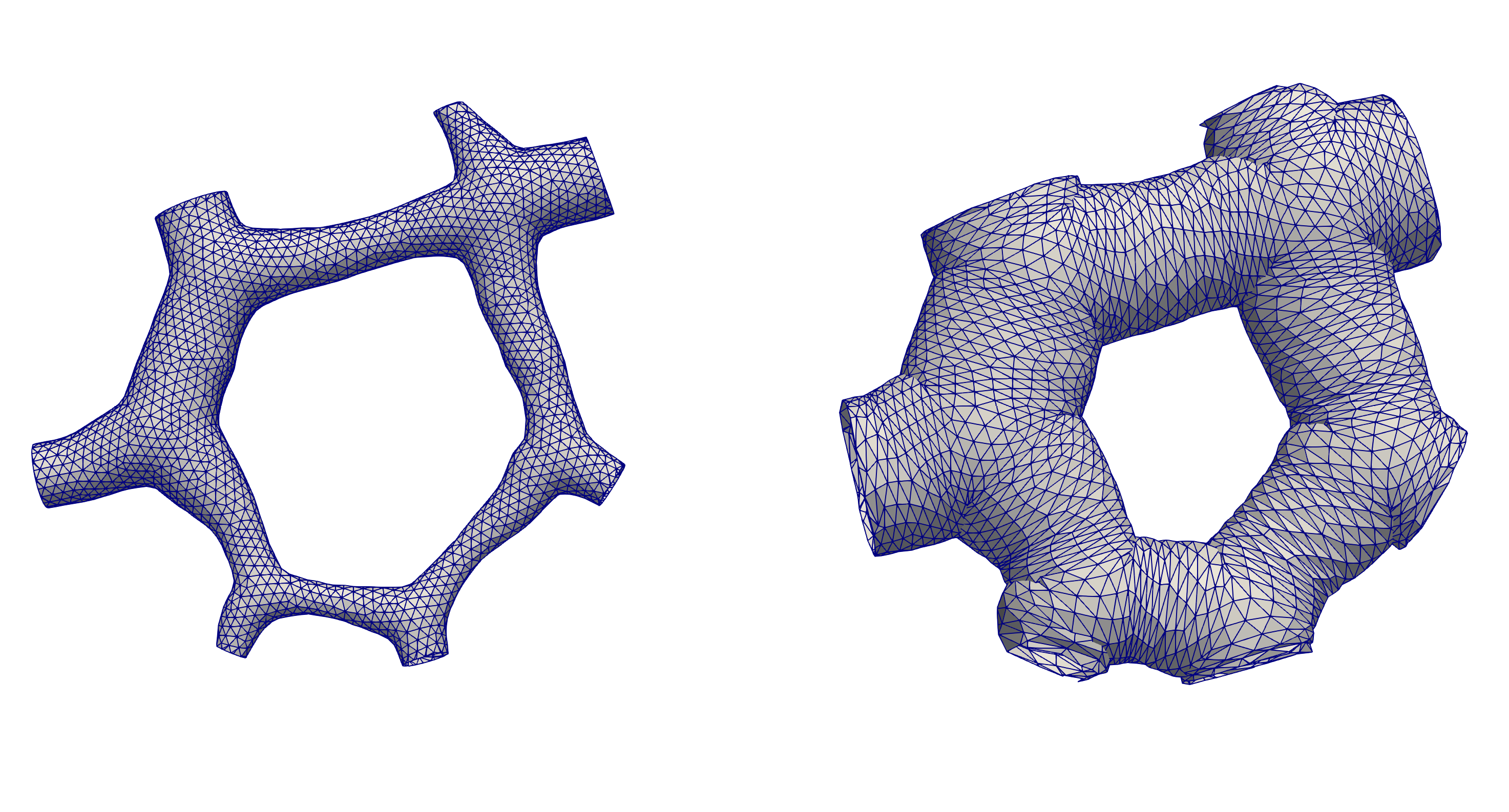}
\put(53,43){\includegraphics[trim = {50cm 3cm 2cm 3cm} , width=0.55\textwidth,  clip=true]{Fig10A.png}}
\put(63,-21){\includegraphics[width=0.162\textwidth]{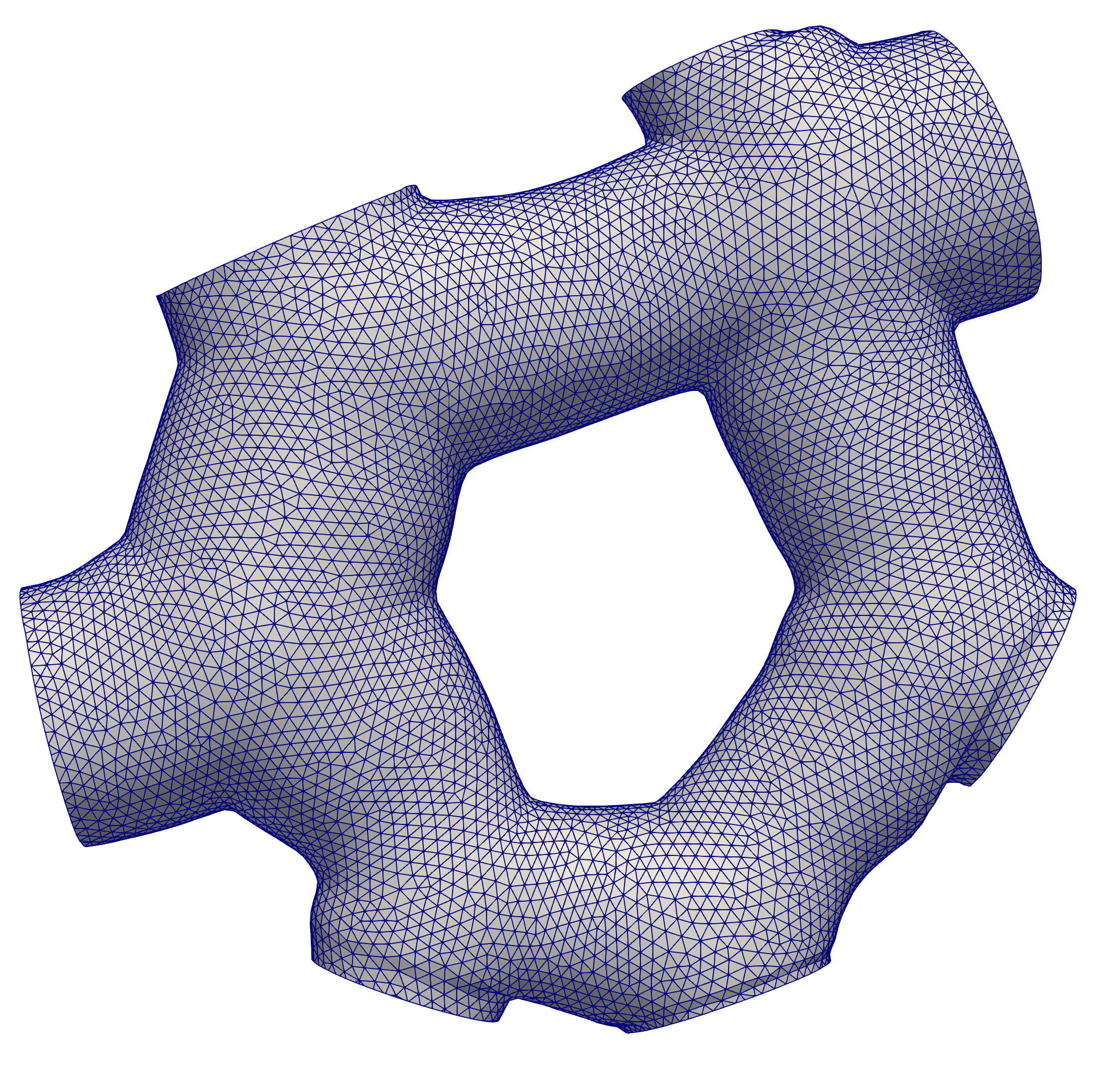}}
\thicklines
\put(49,35){\vector(4,3){14}}
\put(49,26){\vector(4,-3){14}}
\put(44,37){\rotatebox{36}{ Fixed mesh}}
\put(48,30){\rotatebox{-36}{ Remeshed}}
\put(1,90){$a)$}
\put(120,90){$b)$}
\put(120,30){$c)$}
\put(22,-2){$t_0$}
\put(86,-30){$t_1$}
\put(86,39){$t_1$}
\end{overpic}
\end{subfigure}
\begin{subfigure}[t]{0.43\textwidth}
 \vspace{-13cm}
  \vfill
		\includegraphics[trim = {0cm 1cm 4.2cm -10cm}, width=1\textwidth,clip=true, right]{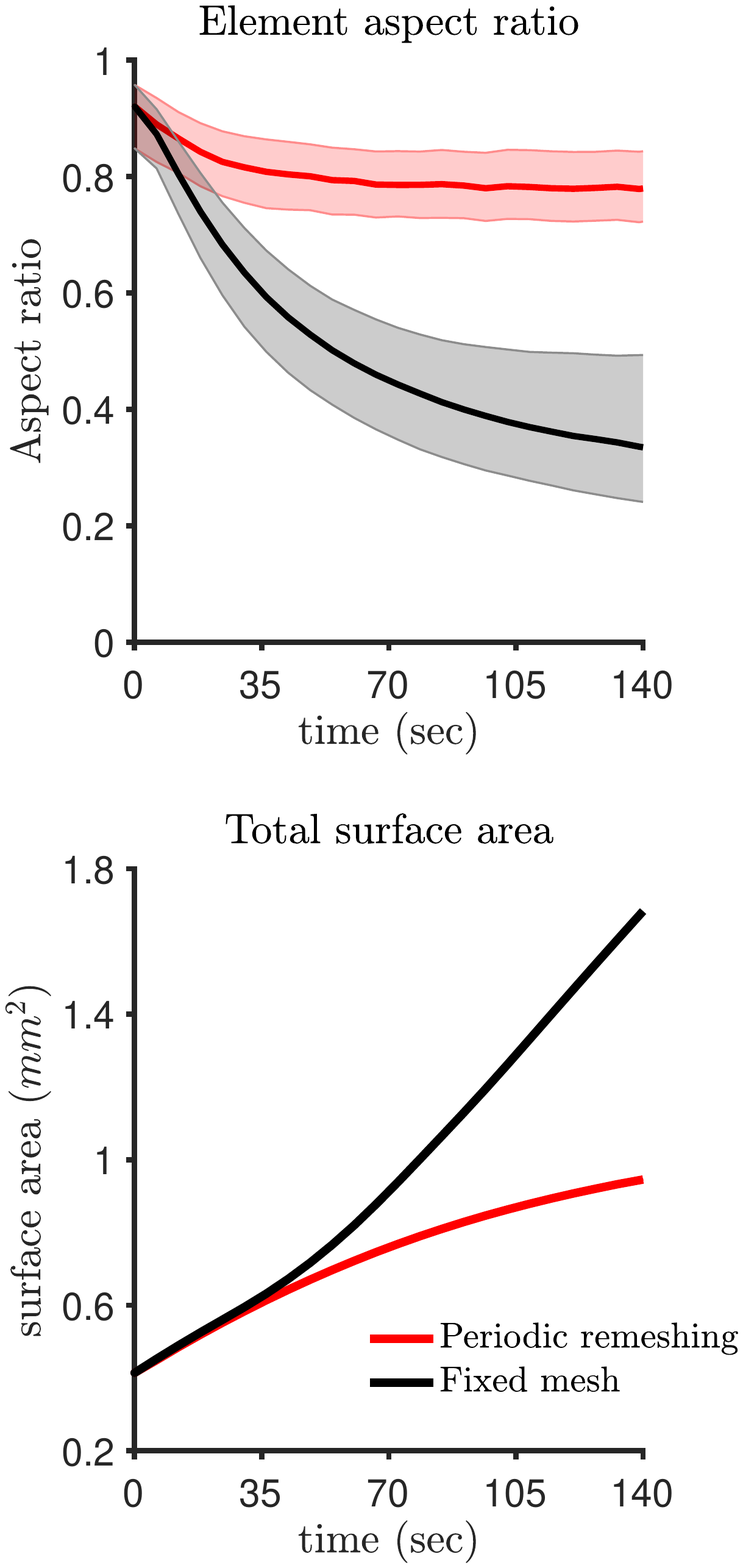}
\end{subfigure}
\caption{\setstretch{1.4}\textbf{A deforming capillary plexus modelled as a hyperelastic membrane with a fixed mesh (upper) or with periodic remeshing (lower).} $(a)$ Periodic remeshing (remeshing every $0.6$ seconds) stabilises the mesh throughout the simulation, whilst the fixed mesh distorts. $(b)$ This is confirmed by the persistently high median ($\pm$ IQR) element aspect ratio with period remeshing (red) and the fall of the element aspect ratio over the fixed mesh (black). 
$(c)$ Without remeshing (black), the total surface area diverges, however with periodic remeshing (red), the total surface area of the plexus grows to an equilibrium.} \label{fig:RemeshingExample}
\end{figure}

\subsection{History dependent remeshing stabilises a hyperelastic model of capillary deformation}
We now turn our attention to the capillary plexus shown in Figure 1, where instability arose from element distortion as the plexus deformed, corrupting the simulation. We model the plexus as a hyperelastic membrane, as described in Section 3.2, deforming with a constant internal pressure, and with either periodic remeshing (remeshing every $0.6$ seconds), or no remeshing (a fixed mesh). Figure \ref{fig:RemeshingExample}a shows the plexuses at time $t_1=140$ seconds both with (lower) and without (upper) remeshing for comparison. It is again evident that without remeshing the mesh undergoes excessive distortion, with the median ($\pm$ IQR) element aspect ratio falling below $0.4$ (Figure \ref{fig:RemeshingExample}b). The plexus modelled with remeshing deforms smoothly, and maintains mesh stability (median element aspect ratio $\sim 0.8$).  The distortion of the fixed mesh introduces error into the simulation, and drives divergence of the plexus surface area (Figure \ref{fig:RemeshingExample}c, black), where convergence is expected, and seen during plexus deformation with remeshing (Figure \ref{fig:RemeshingExample}c, red). 

\section{Discussion} \label{Discussion}
In this paper we have presented a relatively simple history-dependent remeshing technique to stabilise finite element models of surface mesh deformation by removing the error arising with mesh degradation. 
Over a finite element mesh, we consider the domain as a patchwork of continuous linear elements, and utilise this concept to map any point back to its initial position. This mapping enables us to transfer the initial configuration to the new mesh. With this information, we can reconstruct the history dependent variables, such as stress and strain, instead of directly transferring state variables, as is tradition \cite{peric1996transfer,zienkiewicz1992superconvergent, zienkiewicz1992superconvergent2, zienkiewicz1992superconvergent3}. 
These traditional variable transfer remeshing techniques currently available in the literature are highly technical and require significant code implementation \cite{peric1996transfer,zienkiewicz1992superconvergent, zienkiewicz1992superconvergent2, zienkiewicz1992superconvergent3}. The remeshing technique described here provides a simple alternative to these more complex methods, whilst reducing the error accumulated during mesh distortion, and preventing the introduction of degradation related artefacts.  \\

 The mesh configuration is a key consideration when developing and implementing finite element models of deformation. If not handled appropriately, mesh distortion can introduce error into the simulation and sabotage the final results, as shown in Figure \ref{fig:ExampleOfMeshDegredation}.  This matter of mesh distortion is particularly problematic in large, or unidirectional deformations (Figure \ref{fig:ExampleOfMeshDegredation}). 
In this paper we intentionally chose simple examples to analyse this remeshing technique because, while conceptually simple, they demonstrated that this remeshing technique handles 1) element distortion, 2) uniform and unidirectional extension and compression and 3) large deformation, with ease. 
The deformation shown in Figures \ref{fig:2DErrorAnal}--\ref{fig:StrainErrorFrequ} was known, indeed prescribed for the purposes of validation, however it is not normally known (Figure \ref{fig:RemeshingExample}). As such, this technique recovers the reverse deformation for each element and node in simulations where the deformation is not necessarily known. \\

The history-dependent remeshing technique presented here removes mesh distortion related error, however, there is of course some error, although small, inherent to mapping the initial conditions from one mesh to another (Figures \ref{fig:SpatialErrorInDefomredCylinder}--\ref{fig:StrainErrorFrequ}). This error decreases with increasing refinement of the original mesh, $\mathbf{M}$, suggesting that the remeshing error most likely arises from any spatial distance between any given new node in $\widehat{\mathbf{M}}$ and the nearest element from $\mathbf{M}$. Moreover, this error associated with remeshing is sufficiently small compared to that resulting from mesh distortion (Figure \ref{fig:StrainErrorFrequ}a). The remeshing interval is also an important variable influencing the error associated  with remeshing (Figure \ref{fig:StrainErrorFrequ}).  That is to say, the degree of distortion incurred by the original mesh prior to remeshing impacts the accuracy of remeshing. Normally remeshing is triggered by the aspect ratio falling below a given threshold, for example $\sim$ \AR $<0.6$. In most examples shown here, the mesh deformed to three times its original size, which we are considering to be a large deformation. \\

Whilst time consuming when applied to dense meshes, the history dependent remeshing technique described here is more efficient than using an excessive refinement over the initial mesh in an attempt to prevent mesh distortion. Moreover, increasing the mesh refinement alone will be insufficient to prevent element distortion in the case of unidirectional deformation. In addition, the more time consuming step in the remeshing process, that is the element search, could be expedited in future work by the introduction of spatial binning, in which only the elements from within the same predefined spatial bin, or region, as the node are considered. \\

In the development of this history dependent remeshing technique, we have relied on the assumption that all elements deform linearly. This assumption can be relaxed for models assuming higher order deformation across the element, by introducing higher order terms to the basis vectors ($\vec{u}^j_k = [x^j_k(x,y,z),~y^j_k(x,y,z),~z^j_k(x,y,z)]$) to encapsulate the non-linear nature of element deformation. It is expected that this extension will be relatively straight-forward to implement. Finally, this method has been developed and implemented for surfaces in 2D and 3D space, and further consideration should be given to 3D tetrahedral meshes, considering a more specific definition of $u_3$ and utilising the shape functions.\\

\section{Conclusion}
In this paper we have presented a novel history-dependent remeshing technique to prevent mesh distortion in deforming surface meshes. This technique relies on the assumption that the mesh can be considered as a patchwork of continuous elements, and uses this assumption to map the initial configuration from the old distorting mesh, to a new optimised mesh. In doing so, history dependent variables can be recalculated using the new initial conditions in a simple, clean manner. We have shown that the error associated with this technique is small and decreases with increasing refinement of the original mesh, and with an appropriate remeshing frequency. 

\clearpage

\bibliographystyle{ieeetr}
\bibliography{Manuscript.bib}
\end{document}